\documentclass{article}
\usepackage [latin1]{inputenc}
\usepackage{amssymb,amsmath,amsthm}
\usepackage{longtable,booktabs}
\usepackage{graphicx}
\usepackage{amsmath}
\usepackage{mathrsfs}
 \usepackage{cite}
\usepackage{color}
 \usepackage{indentfirst}
\begin{document}
\newcommand{\wzr}[1]{\textcolor{blue}{#1}}
\newtheorem{Lemma}{Lemma}[section]
\newtheorem{Proposition}{Proposition}[section]
\newtheorem{Theorem}{Theorem}[section]
\newtheorem{Corollary}{Corollary}[section]
\newtheorem{Example}{Example}[section]
\baselineskip 0.25in
\title{\Large\bf Hierarchical similarity-based approximate reasoning with restricted equivalence function}
\author{  Dechao Li$^{ \mbox{{\small a}}}$ \qquad Yuhui Zhu$^{ \mbox{{\small b}}}$\thanks{Corresponding Author: Yuhui Zhu; Email: zhuyh\_ld@163.com}\\
 $^{ \mbox{{\small a}}}${\small  \small School of Information and Engineering, Zhejiang Ocean University, Zhoushan, 316000, China}
\\$^{ \mbox{{\small b}}}${\small Shanghai Technical Institute of Electronics and Information, Shanghai, 201411, China}}
\date{}
	\maketitle
	\begin{center}
	\begin{minipage}{140mm}
		\begin{picture}(1,1)
			\line(1,0){400}
		\end{picture}
\centerline{\bf Abstract} \vskip 3mm {\qquad Given that the restricted equivalence functions (REFs) can serve to measure the similarity of two fuzzy sets,  this motivates the integration of REFs with similarity-based approximate reasoning systems to enhance inference capabilities. Therefore, this work primarily constructs hierarchical similarity-based approximate reasoning (SBAR) using REFs. Specifically, we first characterize REFs with a given aggregation function, then discuss the approximation equality of SBAR method proposed by Raha et al. with REFs. Finally, we suggest two REF-based hierarchical Raha's SBAR methods which efficiently restrain the explosion of fuzzy rules.}
 \vskip 2mm\noindent{\bf Key words}: Restricted equivalence function; Fuzzy implication; Aggregation function; Approximation equality; Hierarchical fuzzy reasoning

\begin{picture}(1,1)
 \line(1,0){400}
\end{picture}
\end{minipage}
\end{center}
\vskip 4mm
\section{Introduction}
\qquad The precise distinction or identification of objects constitutes a significant research subject in artificial intelligence. However, practical applications frequently necessitate distinguishing between fuzzy or uncertain objects due to the complexity of human perception and environmental factors. Consequently, distinguishing fuzzy objects by emulating human cognitive processes becomes indispensable across many applied domains such as artificial intelligence. As an extension of metrics, fuzzy metrics was presented to discriminate between two fuzzy sets by Kaleva and Seikkala\cite{Kaleva}. However,  it remains uncertain whether human behavior fully adheres to the triangle inequality\cite{Yearsley}. This motivates the exploration of alternative ideas for comparing two fuzzy sets. Unlike the fuzzy metrics proposed by Kaleva and Seikkala, Liu introduced a distance measure for fuzzy sets\cite{Liu}. Additionally, Liu defined a similarity measure of fuzzy sets and demonstrated the duality of distance measure and similarity measure\cite{Liu}. After, many  similarity measures for fuzzy sets were presented to discriminate two fuzzy sets\cite{Beg, Bustince1,Chen, Baets, Fan, LiQ, Pappis, WangB}. It should be noted that certain similarity measures in the aforementioned literature fail to fulfill the axiomatic definition of similarity measures proposed by Liu. Today, despite widespread application in many fields, selecting an appropriate similarity measure for fuzzy sets remains practically challenging.
\subsection{Motivation}
\qquad As previously stated, Bustince et al. defined a similarity measure  $S(A,B)=M_n(F(A(x_1),$ $B(x_1)), \cdots, F(A(x_n),B(x_n))$ to distinguish two images, where $F$  is a REF and $M_n$ an $n$-ary aggregation function\cite{Bustince1}. They further verified that $S$ satisfies the axiomatic definition of similarity measure proposed by Liu. Obviously, REF acts as a key in the similarity measure of  Bustince et al. We know that  REF strengthens the conditions of Fodor and Roubens' equivalence to meet image comparison requirements\cite{Fodor}. Consequently, REFs constitute a subclass of Fodor-Roubens' equivalences\cite{Bustince1}. As Fodor and Roubens' equivalence of $E$ has been characterized by a fuzzy implication $I$ satisfying IP as $E(x,y)=I(x,y)\wedge I(y,x)$\cite{Fodor}, characterizing REF via fuzzy implications facilitates its application in image comparison.  Bustince et al. similarly characterized  REF  through a fuzzy implication which fulfills OP, CP$(N)$ and CC (see Theorem 7 in \cite{Bustince1}). They also verified that the function $F(x, y)=\varphi^{-1}(1-|\varphi(x)-\varphi(y)|)$ becomes a REF which generated by the fuzzy implication $(I_L)_\varphi(x,y)=\varphi^{-1}((1-|\varphi(x)+\varphi(y)|)\wedge1)$ (see Theorem 8 in \cite{Bustince1}). It is not difficult to find that this REF can also be used to construct the similarity measure $S_L$ proposed by Pappis and Karacapilidis\cite{ Pappis}.

 Recently, Bustince et al. posed an open problem whether the function $F(x,y)=T(I(x,y),$ $I(y,x))$ becomes a REF for a given t-norm $T$ and a fuzzy
implication $I$\cite{Bustince4}. Subsequent research has confirmed that the construction of REFs is achievable not only via t-norms but also through overlap functions\cite{ Palmeira,Qiao}. Obviously, the t-norm and overlap function are two special class of binary aggregation functions\cite{Grabisch}. As fundamental tools in information fusion and decision-making, aggregation functions also serve critical roles in fuzzy logic. Therefore, it is natural extension to investigate  whether the function $F(x,y)=M(I(x,y),I(y,x))$ is a REF for a given binary aggregation function $M$.

The investigation of human reasoning mechanisms has become a prominent research direction of artificial intelligence. As approximate reasoning can effectively simulate human thinking, it has been successfully utilized in data mining, decision-making and artificial intelligence. This success motivates the ongoing pursuit of more suitable approximate reasoning methods. Using the similarity measure, Turksen and Zhong developed Approximate Analogical Reasoning Scheme (AARS)\cite{Turksen}. In this method, the output $B'$ is computed as $B'(y)=m(S(A',A), B(y))$, where $m$ is referred to as a modification function.  Notably, the if-then rules were excluded from the modification function's formulation in AARS. After, Raha et al. introduced SBAR method  which  incorporates the if-then rules as fuzzy relations\cite{Raha}. Obviously, Raha's SBAR method constitutes an extension of Turksen and Zhong's approach. The incorporation of fuzzy relations enhances the effectiveness of Raha's SBAR method. Considering that REF can be employed to compare the  membership degrees of two fuzzy sets, it is feasible to define approximate equality between fuzzy sets using REFs. Consequently, investigating the approximate equality of Raha's SBAR method under REFs is of significant interest.

In practical applications, the rule explosion phenomenon in multi-input-single-output (MISO) fuzzy systems undermines computational efficiency of approximate reasoning. As the sup-projection operation is involved in Raha's SBAR method, computational complexity increases exponentially with the number of fuzzy rules in MISO fuzzy systems employed Raha's SBAR method\cite{Cornelis}. To effectively overcome the rule explosion, some hierarchical approximate reasoning methods have been extensively investigated. For example,  Jayaram proposed a
hierarchical compositional rule of inference (CRI) method\cite{Jayaram}. By leveraging the functional equation $I(T(x,y),z)=I(x,I(y,z))$, an MISO fuzzy system can be transformed into an SISO hierarchical fuzzy system\cite{Jayaram}. Afterward, other hierarchical approximate reasoning methods have been studied\cite{LiG,LiL}. Obviously, the functional equations play crucial roles in these hierarchical approximate reasoning methods. Thus, by investigating the functional equation $S_F(T(A'_1,A'_2),T(A_1,A_2))=T(S_F(A'_1,A_1),S_F(A'_2,A_2))$, it is viable to construct some hierarchical Raha's SBAR methods, where $S_F$ is a similarity measure generated by a REF (detailed in Section 4).
\subsection{Review of related results}
\qquad To compare two images, Bustince et al. introduced the notion of REF and a similarity measure that aggregates REF in 2006\cite{Bustince1,Bustince2}. After, they defined a restricted dissimilarity function (RDF) and an $E_N$-function, then explored the connections between these functions\cite{Bustince3}. In 2010, Chaira introduced an intuitionistic fuzzy REF to compare segmentation of medical images\cite{Chaira}. To develop the construction method of REF, Palmeira et al. characterized it using a commutative aggregation function with a neutral element 1\cite{Palmeira}. To compare radial data, Marco-Detchart et al. defined a restricted radial equivalence function\cite{Detchart}. In 2018,  Palmeira et al. proposed a REF to compare the membership degree of two lattice-valued fuzzy sets\cite{PalmeiraB}. In 2019, Qiao derived a REF from a fuzzy implication and an overlap function\cite{Qiao}. That same year, Altalhi et al. use a REF to construct the similarity measure of data\cite{Altalhi}. In 2020, Bustince et al. proposed an interval-valued REF\cite{BustinceD}. Recently,  Miguel et al. investigated a type-2 REF\cite{Miguel}, while Jaurrieta et al. developed an $n$-dimensional REF\cite{Jaurrieta}. Moreover, Bustince et al. constructed a REF from some strong negations\cite{Bustince4}. These studies are summarized in Table 1.

As a vital concept of fuzzy set theory, the similarity measure has been used to construct some strategies of approximate reasoning\cite{Deng,Raha,Turksen,WangM}. Moreover, these similarity measures were also applied to estimate the performance of approximate reasoning. For example, Li et al. investigated the approximation property properties of Raha's SBAR method using the similarity measure\cite{LiK}. The robustness of some approximate reasoning methods has been investigated  through the similarity measures of fuzzy sets\cite{Dai,Luo,Qin,WangD}. Wang et al. discussed the approximation equality of approximate reasoning\cite{WangB}. However, no existing work has employed REFs to either construct new approximate reasoning strategies or study the approximation equality of Raha's SBAR method.
 $$\mbox{\bf{\small Table\ 1 \ The related results of REF}}$$\vspace{-8mm}
\begin{center}
\begin{tabular}{cc}
 \toprule[1pt]
References &Main contributions\vspace{1mm}\\
\midrule[0.75pt]
Bustince et al.\cite{Bustince1}, 2006& Defined a REF\vspace{1mm}\\
Bustince et al.\cite{Bustince3}, 2008& Studied the relationship among RDF, REF and $E_N$-function\vspace{1mm}\\
Chaira\cite{Chaira}, 2010& Introduced an intuitionistic fuzzy REF\vspace{1mm}\\
Palmeira et al.\cite{Palmeira}, 2013& Characterized the REF with a pseudo t-norm\vspace{1mm}\\
Palmeira et al.\cite{PalmeiraB}, 2018& Investigated the lattice-valued REF\vspace{1mm}\\
Qiao\cite{Qiao}, 2019 &Constructed the REF from overlap function\vspace{1mm}\\
Altalhi et al.\cite{Altalhi}, 2019&Discussed the relationship between deviation function and REF\vspace{1mm}\\
Bustince et al.\cite{BustinceD}, 2020& Proposed an interval-valued REF with admissible orders\vspace{1mm}\\
 Miguel et al.\cite{Miguel}, 2022& Considered a type-2 REF \vspace{1mm}\\
Bustince et al.\cite{Bustince4}, 2022& Constructed the REF from some strong negations\vspace{1mm}\\
 Jaurrieta et al.\cite{Jaurrieta}, 2023& Defined a REF on $[0,1]^n$\vspace{1mm}\\
   \bottomrule[1pt]
\end{tabular}
\end{center}\vspace{1mm}
\subsection{Contribution of this work}
\qquad As discussed above, REFs can be used to construct the similarity measure of two fuzzy sets. To enhance the applicability Raha's SBAR method in fields like artificial intelligence, image processing and classification, this work will integrate REF with Raha's SBAR method, and study the performance of Raha's SBAR method. The primary contributions of this work are centered on

(1) Characterizing REFs with a given binary aggregation function.

(2) Investigating the approximation equality of Raha's SBAR method with REFs.

(3) Constructing two hierarchical Raha's SBAR methods with REFs.

Having these in mind, we arrange this work as follows. Section 2 reviews some fundamental concepts such as fuzzy negation, aggregation function,  fuzzy implication, REF and Raha's SBAR method. Section 3 characterizes REFs for a given aggregation function. In Section 4, the approximation equality of Raha's SBAR method is studied using REFs. Section 5 presents two hierarchical Raha's SBAR methods with REFs.
\section{Preliminary}
\qquad In this section, we briefly summarize some basic concepts and results needed for further treatment. Firstly, let us recall the negation, aggregation function and fuzzy implications.
\subsection{Negation, aggregation function and fuzzy implication}
{\bf Definition 2.1}\cite{Lowen} A non-increasing mapping $N:[0,1]\rightarrow [0,1]$ is called
a fuzzy negation if  $N(0)=1$ and $N(1)=0$. Further,

i. $N$ is strict when it is continuous and strictly decreasing;

ii. $N$ is strong when $N(N(x))=x$ for any $x\in [0,1]$.
\\{\bf Example 2.2}[17] The standard negation $N_0(x)=1-x$  is a strong fuzzy negation.
\\{\bf Definition 2.3}\cite{Grabisch} The mapping $M_n:[0,1]^n\rightarrow [0,1]$ is called an $n$-ary aggregation function if

(A1) $M_n(0, 0,\cdots, 0)=0$ and $M_n(1, 1,\cdots, 1)= 1$,

(A2)  $M_n(x_1, x_2,\cdots, x_n)\leq M_n(y_1, y_2,\cdots, y_n)$ if $x_i\leq y_i(i=1,2,\cdots,n)$.
\\{\bf Definition 2.4}\cite{Grabisch} Let $M_n$ be an $n$-ary aggregation function.

i.  $a$ is an annihilator of $M_n$ if $M_n(x_1,\cdots,x_{i-1},a,x_{i+1},\cdots,x_n)=a$ holds for each $n\geq 2$, each $i\in \{1,2,\cdots,n\}$ and for all $x_1,\cdots,x_{i-1},x_{i+1},\cdots,x_n\in[0,1]$;

ii. $e\in[0,1]$ is a neutral element of $M_n$ if we have $M_n(x_1,\cdots,x_{i-1},e,x_{i+1},\cdots,x_n)=M_n(x_1,\cdots,x_{i-1},x_{i+1},\cdots,x_n)$ for each $n\geq 2$, each $i\in \{1,2,\cdots,n\}$ and for all $x_1,\cdots,x_{i-1},x_{i+1},$ $\cdots,x_n\in[0,1]$.
\\{\bf Definition 2.5}\cite{Grabisch} The $n$-ary aggregation function $M_n$ is said to be

i. symmetric if $M_n(x_1,x_2,\cdots, x_n) = M_n(x_{p(1)},x_{p(2)},\cdots,x_{p(n)})$ for any permutation $p$ on $\{1,2\cdots,n\}$;

ii. conjunctive if $M_n(x_1,x_2,\cdots, x_n)\leq x_1\wedge x_2\wedge \cdots\wedge x_n$;

iii. disjunctive if $M_n(x_1,x_2,\cdots, x_n)\geq x_1\vee x_2\vee \cdots\vee x_n$;

iv. averaging if $x_1\wedge x_2\wedge \cdots\wedge x_n\leq M_n(x_1,x_2,\cdots, x_n)\leq x_1\vee x_2\vee \cdots\vee x_n$.
\\{\bf Remark 1.} For convenience to write, the binary aggregation function $M_2$ will be abbreviated as $M$ in the rest of this paper.
\\{\bf Definition 2.6}\cite{Klement} A binary aggregation function $T$ is called a t-norm if it is commutative, associative, non-decreasing and has a neutral element 1.
\\{\bf Example 2.7}\cite{Klement} The four basic t-norms are respectively

$T_M(x,y)=x\wedge y$,

$T_P(x,y)=xy$,

 $T_L(x,y)=(x+y-1)\vee0$,\vspace{2mm}

$T_D(x,y)=\left\{\begin{array}{ll}
             x\wedge y& x=1\ \textmd{or}\ y=1\\
             0& \textmd{otherwise}
           \end{array}\right.$.\vspace{2mm}
\\{\bf Definition 2.8}\cite{Klement} i.  $x\in (0,1)$ is called a nilpotent element of $T$ if  $x_T^{(n)}=0$ for some $n\in\mathbb{N}$, where $x_T^{(n)}=T(x,x,\cdots,x)$.

ii. $x\in (0,1)$ is called a zero divisor of $T$ if there exists some $y\in(0,1)$ such that  $T(x,y)=0$.
\\{\bf Definition 2.9}\cite{Klement} i. $T$ is called Archimedean if an $n\in \mathbb{N}$ can be found such that $x^{(n)}_T<y$ for all $x$ and $y\in (0,1)$;

ii. $T$ is called nilpotent if it is continuous  and each $x\in (0,1)$ is its nilpotent element;

iii. $T$ is called strict if it is continuous and strictly monotone (that is, $T(x,y)<T(x,z)$ holds for all $x>0$ and $y<z$).
\\{\bf Proposition 2.10}\cite{Klement} i. $T$ is strict $\Longleftrightarrow$  $T=(T_P)_t$, where $(T_P)_t(x,y)=t^{-1}(t(x)t(y))$ and $t$ is  an increasing bijection on $[0,1]$;

ii. $T$ is nilpotent $\Longleftrightarrow$  $T=(T_L)_t$.
\\{\bf Definition 2.11}\cite{Baczynski,Dimuro,Palmeira} A mapping $I:[0,1]^2\rightarrow [0,1]$ is called a fuzzy implication if it fulfills

(I1) non-increasing in the first variable,

(I2) non-decreasing in the second variable,

(I3) $ I(0,0)=I(0, 1)=1,\ I(1, 0)=0$.

Obviously, $I(0,y) = 1$ and $I(x,1)=1$ hold for any $x, y\in[0, 1]$. Further, we say that $I$ satisfies

(CC) $I(x,y) = 0\Longleftrightarrow x = 1, y = 0$;

(CP($N$)) $I(x, y) = I(N(y),N(x))$;

(EP)  $I(x, I(y, z)) = I(y, I(x, z))$;

(IP) $I(x, x) = 1$;

(NP) $I(1, y) = y$;

(OP) $I(x, y)=1\Longleftrightarrow x\leq y$ for all $x, y\in[0, 1]$;

(LOP) $ x\leq y \Longrightarrow I(x, y)=1$.
 \\{\bf Definition 2.12}\cite{Baczynski} An R-implication $I_T$ associated with  $T$  is
defined by $I_T(x, y)= \textmd{sup}\{a\in[0,1] \mid T(x, a)\leq
y\}$.
For example, the Goguen implication  $I_{GG}(x,y)=\frac{x}{y}\wedge 1$ is generated by $T_P$,  the {\L}ukasiewicz implication  $I_{L}(x,y)=(1-x+y)\wedge 1$ is generated by $T_P$ and the G\"{o}del\vspace{1mm} implication $I_G(x,y)=\left \{\begin{array}{lc}
  1,&  x\leq y
\\y, & x>y
\end{array}\right.$ is generated by $T_M$.\vspace{1mm}
\\{\bf Proposition 2.13}\cite{Baczynski} i. $I_T$ fulfills NP and IP. Further, $I_T$ satisfies EP, OP and is right continuous for the second variable if $T$ is left continuous;

ii. $T(x,y)\leq z\Longleftrightarrow x\leq I_T(y,z)$ if $T$ is left continuous.
\subsection{REF, similarity measure and SBAR}
\qquad This subsection will introduce the concept of REF and similarity measure, while also recalling Raha's SBAR method.
\\{\bf Definition 2.14}\cite{Bustince1}  The mapping $F:[0, 1]^2\rightarrow [0,1]$ is referred as a restricted equivalence function (REF) if

(REF1) $F(x, y) = F(y, x)$,

(REF2) $F(x, y) = 1\Longleftrightarrow x = y$,

(REF3) $F(x, y) = 0 \Longleftrightarrow x = 1, y = 0$ or $x = 0, y = 1$,

(REF4) $F(x, y) = F(N(x), N(y))$ with a strong negation $N$,

(REF5) $F(x, z)\leq F(x, y)$ if $x\leq y\leq z$.

 Let $\mathcal{F}(U)$ denote the family of fuzzy sets  on the universe $U$.
\\{\bf Definition 2.15}\cite{Liu, Raha}  $S:\mathcal{F}(U)\times \mathcal{F}(U)\rightarrow [0, 1]$ is called a similarity measure if

(S1) $S(A, B)=S(B, A)$,

(S2) $S(A, B)=1\Longleftrightarrow A=B$,

(S3) $S(A,B)=0$ implies $A(x)\wedge B(x)=0$ for all $x\in U$,

(S4) $A\subseteq B\subseteq C$ implies $S(A, C)\leq S(A, B)\wedge S(B,C)$.

Finally, we recall the SBAR method proposed by Raha et al.\cite{Raha}. In the Raha's SBAR method,  four steps are given to calculate the conclusion of fuzzy reasoning. Practically speaking,  the conditional statement firstly is interpreted  as a conditional fuzzy relation $R(A, B)$. Secondly, the similarity degree $S(A, A')$ between the input $A'$ and the antecedent of a fuzzy rule $A$ is calculated. Further,  $S(A, A')$ is modified as a conditional relation $R(A, B|A')$. Finally, the conclusion is computed by using the sup-projection operation on $R(A, B|A')$. For implication-based and conjunction-based of if-then rules, the conclusions $B'_1$ and $B'_2$ can be respectively expressed as follows
\begin{equation}B'_1(y)=\mathop{\sup}\limits_{x\in U}I_T(S(A',A),T(A(x),B(y))),
\end{equation}
\begin{equation}B'_2(y)=\mathop{\inf}\limits_{x\in U}I_T(S(A',A),I_T(A(x),B(y))),
\end{equation}
where $I_T$ is an R-implication and $S$ is a similarity measure.
\section{Characterization of REFs with a given binary aggregation function}
\qquad We know that REFs can be constructed from t-norms and fuzzy implications. By generalizing t-norms to binary aggregation functions in this section,  this section will investigate whether REFs can be constructed using a binary aggregation function.

Let $M$ denote a given binary aggregation function, i.e., a binary operation on the interval $[0,1]$ fulfilling A1 and A2. Using $M$, we can define a mapping $G:[0,1]^2\rightarrow [0,1]$ as $G(x,y)=M(f(x,y),f(y,x))$ for all $x,y\in[0,1]$, where $f:[0,1]^2\rightarrow [0,1]$ is an arbitrary mapping. We subsequently determine the additional requirements needed for $f$ and $M$  to ensure that $G$ fulfills properties REF1 through REF5 respectively.
\\{\bf Proposition 3.1} $G$ fulfills REF1 if $M$ is commutative. Conversely, $M$ is commutative if $G$ satisfies REF1 and $f$ is surjective and continuous on $[0,1]^2$.
\\{\bf Proof.} Obviously.
\\{\bf Proposition 3.2} Let $M$ be one strict, that is, $M(x,y)=1\Longrightarrow x=1$ and $y=1$. $G$ fulfills REF2 if $f$ satisfies OP.
\\{\bf Proof.}  Suppose that $f$ satisfies OP. We have $f(x,y)=1\Longleftrightarrow x\leq y$. This implies that $G(x,y)=M(f(x,y),f(y,x))=M(1,1)=1$ if $x=y$. Conversely, we assume  that $G(x,y)=M(f(x,y),f(y,x))=1$. Since $M$ is one strict, $f(x,y)=f(y,x)=1$ hold. Thus, we have $x\leq y$ and $y\leq x$. That is, $x=y$.
\\{\bf Remark 2.} Obviously, $M$ is one strict if it has a neutral element 1. Therefore, Proposition 3.2 holds for all binary aggregation functions with a neutral element 1.
\\{\bf Proposition 3.3} Let $M$ be an averaging  aggregation and $f$ satisfy I1. $f$ fulfills OP if $G$ fulfills REF2.
\\{\bf Proof.} To ensure the fact that $f$ satisfies OP, we first verify $f(x,y)=1$ when $x\leq y$. As $G$ fulfills REF2, we have $G(y,y)=M(f(y,y),f(y,y))=1$. This means that $f(y,y)=f(y,y)\vee f(y,y)\geq M(f(y,y),f(y,y))=1$ holds. Thus, $I(x,y)\geq I(y,y)=1$. Next, we confirm that either $x\leq y$ if $f(x,y)=1$ or  $y\leq x$ if $f(y,x)=1$. Otherwise, we assume that $x\neq y$. Since $M$ is averaging, $f(x,y)\wedge f(y,x)\leq M(f(x,y),f(y,x))<1$ holds. This implies that $f(x,y)<1$ or $f(y,x)<1$. This is a contradiction.
\\{\bf Proposition 3.4} Let $f$ satisfy I1 and $f_M(y)=M(1,y)$ be strictly increasing on [0,1]. $f$ satisfies OP if $G$ fulfills REF2.
\\{\bf Proof.} Similarly, we can verify $f(x,y)=1$ if $x\leq y$. Thus, it is enough to ensure that $x\leq y $ when $f(x,y)=1$. Otherwise, suppose that $x>y$.
This case implies that $G(x,y)=M(f(x,y),f(y,x))=M(1,f(y,x))=f(y,x)\geq f(x,x)=1$. The REF2 implies that $x=y$. This is a contradiction.
\\{\bf Proposition 3.5}  Let $M$ be a binary aggregation function without zero divisors and $M(0,1)=M(1,0)=0$. $G$ satisfies REF3 if $f(0,1)=1$ and $f$ fulfills CC.
\\{\bf Proof.} Let us firstly verify that $G(0,1)=0$ or $G(1,0)=0$. Since $f(0,1)=1$ and  $f$ fulfills CC,  $G(1,0)=M(f(1,0),f(0,$ $1))=M(0,1)=0$ holds. Similarly, we have $G(0,1)=0$.
Inversely, suppose that $G(x,y)=0$. We then have $G(x,y)=M(f(x,y),f(y,x))=0$. The fact that $M$ has not zero divisors implies that $f(x,y)=0$ or $f(y,x)=0$ holds.  Therefore, we have $x=1,y=0$ or $x=0,y=1$ according to CC. Thus, $G$ satisfies REF3.
\\{\bf Proposition 3.6} Let 1 be a neutral element of $M$.  $G$ satisfies REF3 if $f$ fulfills CC and LOP.
\\{\bf Proof.} Let us firstly verify that $G(0,1)=0$ or $G(1,0)=0$. As 1 is a neutral element of $M$, we have $G(1,0)=M(f(1,0),f(0,1))=M(0,1)=0$. Similarly, $G(0,1)=0$ holds.
Conversely, suppose that $G(x,y)=M(f(x,y),f(y,x))=0$. Let us consider the following two cases:
i. $x\leq y$. We then obtain $0=G(x,y)=M(f(x,y),f(y,x))=M(1,f(y,x))=f(y,x)$. Thus, $x=0$ and $y=1$ hold by CC.
ii. $x\geq y$. This implies $0=G(x,y)=M(f(x,y),f(y,x))=M(f(x,y),1)=f(x,y)$. We then have $x=1$ and $y=0$  by CC. Therefore,  $G$ satisfies REF3.

The following statements establish the conditions for $f$ to satisfy CC when $G$ fulfills REF3.
\\{\bf Proposition 3.7} Let $M$ be a binary aggregation function without zero divisors and satisfy $M(0,1)=0$. Then, $f$ satisfies CC if $G$ fulfills REF3 and $f(0,1)\neq 0$.
 \\{\bf Proof.}  Let us first verify that $f(1,0)=0$. According to REF3, we have $G(1,0)=M(f(1,0),f(0,$ $1))=0$  or $G(0,1)=M(f(0,1),f(1,0))=0$. As $M$ have no zero divisor,  $f(1,0)=0$ or $f(0,1)$ holds. Thus, $f(1,0)=0$.
On the other hand, suppose that $f(x,y)=0$. This implies $G(x,y)=M(f(x,y),f(y,x))= M(0,f(y,x))\leq M(0,1)=0$. By REF3, we have $x=1,y=0$ or $x=0,y=1$. $f(0,1)\neq 0$ means $x=1,y=0$. Therefore, $f$ fulfills CC.
\\{\bf Proposition 3.8}  Let 1 be a neutral element of $M$ and $M$ have no zero divisors. $f$ satisfies CC if $G$ fulfills REF3 and $f(0,1)\neq 0$.
\\{\bf Proof.} This proof can be obtained similarly to that of Proposition 3.7.
\\{\bf Proposition 3.9}  $G$ satisfies REF4 if $f$ fulfills CP($N$).
\\{\bf Proof.} Since $f$ satisfies CP($N$), we have $G(N(y),N(x))=M(f(N(y),N(x)),f(N(x),N(y)))$ $=M(f(x,y),f(y,x))=G(x,y)$.
The following proposition shows a condition under which $f$ satisfies CP($N$) when $G$ fulfills REF4.
\\{\bf Proposition 3.10} Let $f_M(y)=M(1,y)$ and $g_M(x)=M(x,1)$ be two strictly increasing mappings from [0,1] to [0,1] and $f$ fulfill LOP.  $f$ fulfills CP($N$) if $G$ satisfies REF4.
\\{\bf Proof.} Let $G$ fulfill REF4. And then we consider the following two cases.
i. $x\leq y$. Since $f$ satisfies LOP,  $f(x,y)=1$ holds. Therefore, we have $G(N(y),N(x))=1$ in this case. This implies that  $G(x,y)=M(f(x,y),f(y,x))=M(1,f(y,x))=f_M(f(y,x))=G(N(y),N(x))=M(f(N(y),N(x)),f(N(x),N(y)))=M(1,f(N(x),N(y)))=f_M(f(N(x),$ $N(y))$. Since $f_M(y)=M(1,y)$ is strictly increasing
on [0,1], we can obtain  $f(y,x)=f(N(x),$ $N(y))$.
  ii. $y\leq x$. This means that $N(x)\leq N(y)$ holds. We then have $f(y,x)=1= f(N(x),N(y))$ by LOP.  Therefore, $G(x,y)=M(f(x,y),f(y,x))=M(f(x,y),1)=g_M(f(x,y))=G(N(y),$ $N(x))=M(f(N(y),N(x)),f(N(x),N(y)))=M(f(N(y),N(x)),1)=g_M(f(N(y),N(x))$. As $g_M(x)=M(x,1)$ is strictly increasing
on [0,1], we have $f(x,y)=f(N(y),N(x))$.
Thus,  $f$ satisfies CP($N$).
\\{\bf Remark 3.} It is evident that $M$ with neutral element 1 satisfies the condition of Proposition 3.10. This means that Proposition 3.10 holds for all binary aggregation functions having a neutral element 1.
\\{\bf Proposition 3.11}  $G$ satisfies REF5 if $f$ fulfills I1 and LOP.
\\{\bf Proof.}  Let $x\leq y\leq z$. As $f$ satisfies I1 and LOP, we have $G(x,z)=M(f(x,z),f(z,x))=M(1,f(z,x))\leq M(1,f(y,x))=M(f(x,y),f(y,x))=G(x,y)$. Here, we use the inequality $f(z,x)\leq f(y,x)$.

The following proposition shows a condition for $f$ to fulfill I1 when $G$ satisfies REF5.
\\{\bf Proposition 3.12} Let $f_M(y)=M(1,y)$ is strictly increasing on [0,1] and $f$ satisfy LOP.  $I$ fulfills I1 if $G$ satisfies REF5.
\\{\bf Proof.} Suppose that $x\leq y\leq z$ and $G(x,z)\leq G(x,y)$ hold. We therefore obtain $M(1,f(z,x))=M(f(x,z),f(z,x))=G(x,z)\leq G(x,y)=M(f(x,y),f(y,x))=M(1,f(y,x))$. Since $f_M(y)=M(1,y)$ is strictly increasing, we have $f(z,x)\leq f(y,x)$. Thus, $f$ fulfills I1.
Building upon the arguments above, we can construct a restricted equivalence function $F$ as $F(x,y)=M(f(x,y), f(y,x))$ with a binary aggregation function $M$ and a mapping $f$ as formalized in the subsequent theorems.
\\{\bf Theorem 3.13}  Let $M$ be a commutative and one strict binary aggregation function without zero divisors that fulfills $M(0,1)=0$. If the mapping $f:[0,1]^2\rightarrow [0,1]$ satisfies I1, CC, CP($N$) and OP, then $G(x,y)=M(f(x,y), f(y,x))$ is a REF.
\\{\bf Proof.} By Propositions 3.1, 3.2, 3.5, 3.9 and 3.11, we can obtain the fact that $G(x,y)=M(f(x,y), f(y,x))$ fulfills REF1-REF5. Thus, $G$ is a REF.
\\{\bf Theorem 3.14} Let $M$ be a commutative binary aggregation function with  a neutral element 1. If the mapping $f:[0,1]^2\rightarrow [0,1]$ satisfies I1, CC, CP($N$) and OP, then $G(x,y)=M(f(x,y), f(y,x))$ is a REF.
\\{\bf Proof.} This can be verified similarly to  Theorem 3.13.
\\{\bf Remark 4.} Obviously, $f$ becomes a fuzzy implication when it satisfies I1, CC, CP($N$) and OP\cite{Baczynski}. Therefore,  fuzzy implications can be utilized to construct REFs. Indeed, Qiao constructed a REF by a fuzzy implication\cite{Qiao}. However, the neutral element is not necessary for $M$ in Theorem 3.13. Consequently, the results in \cite{Qiao} constitute a special case of Theorem 3.13.
\\{\bf Example 3.15} Let $M(x,y)=\frac{x+y}{2}$ and $f$ be the {\L}ukasiewicz implication. By Theorem 3.13, we can construct a REF $F_1$ as
$$F_1(x,y)=\frac{(1-x+y)\wedge1+(1-y+x)\wedge 1}{2}.$$

In the final part of this section, we investigate the construction of the mapping $f$ from a REF $F$ and a binary aggregation function $M$ such that $F(x,y)=M(f(x,y),f(y,x))$.
\\{\bf Theorem 3.16} Let $M$ be a commutative binary aggregation function and let $f_M(y)=M(1,y)$ be strictly increasing on [0,1]. If $F$ is a REF, then there exists a mapping $f:[0,1]^2\rightarrow [0,1]$ which fulfills I1, CC, CP($N$) and OP such that $F(x,y)=M(f(x,y), f(y,x))$.
\\{\bf Proof.}  Suppose that $F$ is a REF. Let us define a mapping $f$ as
$$f(x,y)=\left\{\begin{array}{ll}
 1 & x\leq y\\
 f_M^{(-1)}(1, F(x,y))&\textmd{otherwise} \end{array}\right.,$$ where $f_M^{(-1)}$ is the pseudo-inverse of $f_M$.\vspace{1mm}
 Firstly, we can assert that $F(x,y)=M(f(x,y),f(y,x))$. Indeed, we have $F(x,x)=1=M(1,1)=M(f(x,x),f(x,x))$. Thus, we further consider the case when $x\neq y$. Without the loss of generality, let $x<y$. This means that $M(f(x,y),f(y,x))=M(1,f(y,x))=f_M(f(y,x))=f_M( f_M^{(-1)}(1, F(y,x)))$ $=F(y,x)=F(x,y)$ holds. To accomplish the proof, we need to verify that $f$ fulfills I1, CC, CP($N$) and OP. Obviously, $f$ fulfills LOP. Therefore, $f$ satisfies I1 by Proposition 3.12. On the other hand, $f$ fulfills CP($N$) and OP according to Propositions 3.4 and 3.10.
Thus, it is sufficient to ensure that $f$ satisfies CC. Let $f(x,y)=0$. We have  $f_M^{(-1)}(1, F(x,y))=0$. This means that $F(x,y)=f_M(0)=0$. By REF3,  $x=1$ and $y=0$ hold. Conversely, we have $f(1,0)=f_M^{(-1)}(1, F(1,0))=f_M^{(-1)}(1, 0)=0$. Therefore, $f$ fulfills CC.
\\{\bf Remark 5.} i. Obviously, when 1 is a neutral element of $M$, $f_M(y)=M(1,y)=y$ satisfies the\vspace{1mm} condition of Theorem 3.16. In this case, we have \vspace{1mm}$f(x,y)=\left\{\begin{array}{ll}
1 & x\leq y\\
 F(x,y)&\textmd{otherwise} \end{array}\right.$.\vspace{1mm}

ii. The aforementioned case also appeared in \cite{Palmeira}. However, Theorem 3.16 demonstrates that the neutral element condition is not strictly required for the aggregation function $M$. On the other hand, the results presented by Palmeira et al. can be viewed as a particular instance of Theorem 3.16.
\\{\bf Example 3.17} Let $F_2(x,y)=1-|x-y|$ be a REF and $M(x,y)=\frac{x+y}{2}$. According to Theorem 3.16, we can obtain a mapping $f$ as
$$f(x,y)=\left\{\begin{array}{ll}
1 & x\leq y\\
(1-2|x-y|)\vee 0&\textmd{otherwise} \end{array}\right.$$ such that $F_2(x,y)=\frac{f(x,y)+f(y,x)}{2}$.

As shown in Theorem 3.16, the fuzzy implication which fulfills CC, CP($N$) and OP serve as fundamental tools for constructing REFs. However, the CP($N$) property fails to hold for many well-known fuzzy implications  when employing a strong negation $N$. Consequently, we propose a modification of these fuzzy implications to ensure CP($N$) compliance, which will subsequently be used to construct REF.
\\{\bf Definition 3.18}\cite{Aguilo,Baczynski} Let $N$ be a strong negation and $f:[0,1]^2\rightarrow [0,1]$ be a mapping. The mappings $f^u_N$, $f^l_N$, $f^m_N$ and $f^{l_c}_N$ are respectively defined  as
$$f^u_N(x,y)=\max(f(x,y),f(N(y),N(x))), f^m_N(x,y)=\min(f(x,y)\vee N(x),f(N(y),N(x))\vee y),$$
$$f^l_N(x,y)=\min(f(x,y),f(N(y),N(x))),\qquad f^{l_c}_N(x,y)=\left\{\begin{array}{ll}
f(N(y),N(x)) & y<N(x)\\
f(x,y)&\textmd{otherwise} \end{array}\right..\qquad$$
{\bf Proposition 3.19}\cite{Aguilo,Baczynski} $f^u_N$, $f^l_N$, $f^m_N$ and $f^{l_c}_N$ fulfills I1, CP($N$) and OP when $f$ satisfies I1, CP($N$) and OP.
\\{\bf Proposition 3.20} $f^u_N$, $f^l_N$, $f^m_N$ and $f^{l_c}_N$ fulfills CC if $f$ satisfies CC.
\\{\bf Proof.} We only prove that $f^u_N$ satisfies CC. The others can be similarly verified. Let $f^u_N(x,y)=0$. We then have $ f(x,y)=f(N(y),N(x))=0$. Since $f$ fulfills CC, $x=1,y=0$ holds. This means that $f^u_N$ satisfies CC.
Based on Propositions 3.19 and 3.20, Theorems 3.13 and 3.14 can be rewritten as
\\{\bf Theorem 3.21} Let $M$ be a commutative and one strict binary  aggregation function without zero divisors that fulfills $M(0,1)=0$. If the mapping $f:[0,1]^2\rightarrow [0,1]$ satisfies I1, CC and OP, then $F(x,y)=M(f^*_N(x,y), f^*_N(y,x))$ is a REF, where $*\in \{u,m,l,l_c\}$.
\\{\bf Theorem 3.22} Let $M$ be a commutative binary aggregation function with  a neutral element 1. If the mapping $f:[0,1]^2\rightarrow [0,1]$ satisfies I1, CC and OP, then $F(x,y)=M(f^*_N(x,y), f^*_N(y,x))$ is a REF.
\\{\bf Example 3.23} Let $M(x,y)=xy$. With the Goguen implication $I_{GG}$ and  the standard negation $N_0$, we can construct two REFs $F_3$ and $F_4$ by Theorem 3.22 as
$$F_3(x,y)=\left\{\begin{array}{ll}
1 & x=y\vspace{1mm}\\
\frac{x}{y}\wedge\frac{1-y}{1-x}& x<y\vspace{1mm}\\
\frac{y}{x}\wedge\frac{1-x}{1-y}& x>y
 \end{array}\right.,\qquad F_4(x,y)=\left\{\begin{array}{ll}
1 & x=y\vspace{1mm}\\
\frac{1-x}{1-y}& y<x\wedge(1-x)\vspace{1mm}\\
\frac{y}{x}& 1-x\leq y<x\vspace{1mm}\\
\wedge\frac{1-x}{1-y}& x<y\wedge(1-y)\vspace{1mm}\\
 \frac{x}{y}& 1-y\leq x<y\end{array}\right..$$
\section{Approximate equality of SBAR method with REF}
\qquad In this section, we will define the approximate equality of two fuzzy sets using REFs, then will analyze the retainment of approximate equality under fuzzy set operations. Specifically, we will examine REF-based approximate equality for fundamental operations (complement, union, intersection and composition). Finally, we will investigate the approximation equality of Raha's SBAR method with REFs.
\\{\bf Definition 4.1}\cite{WangB} Let $S$ be a similarity measure on $\mathcal{F}(U)$. We say that two fuzzy sets $A$ and $B$ on the universe $U$ are equal to degree $\alpha$ with respect to $S$ if $S(A,B)\geq \alpha$. They then are denoted as $A\sim^S_\alpha B$.

Let $F$ be a REF. As established in Section 3, $F$ can be generated by an aggregation function $M$ and a fuzzy implication $I$, that is, $F$ can be expressed as $F(x,y)=M(I(x,y),I(y,x))$. We define a similarity measure $S_F$ as $S_F(A,B)=\mathop{\inf}\limits_{x\in U}F(A(x),B(x))$ which extends the approach in \cite{Bustince1}. Consequently, the equality to degree $\alpha$ with respect to $S_F$ can be denoted by $A\sim^{S_F}_\alpha B$.. Further, we assume that $M$ and $I$ satisfy the conditions of Theorems 3.13 (3.14) or 3.21 (3.22). This yields the representation $S_F(A,B)=\mathop{\inf}\limits_{x\in U}M(I(A(x),B(x)),I(B(x),A(x)))$. Then, we have the following statements.\vspace{1mm}
\\{\bf Proposition 4.2} Let $A^C$ denote the complement of fuzzy set $A$ and $A^C(x)=N(A(x))$ with a strong negation $N$. We have $A^C\sim^{S_F}_\alpha B^C$ if $A\sim^{S_F}_\alpha B$.
\\{\bf Proof.} Let $A\sim^{S_F}_\alpha B$. We have  $S_F(A^C,B^C)=\mathop{\inf}\limits_{x\in U}F(A^C(x),B^C(x))=\mathop{\inf}\limits_{x\in U}F(N(A(x)),$ $N(B(x)))=\mathop{\inf}\limits_{x\in U}F(A(x),B(x))=S_F(A,B)\geq \alpha$. Thus, $A^C\sim^{S_F}_\alpha B^C$.\vspace{1mm}
\\{\bf Proposition 4.3} $A\cup C\sim^{S_F}_\alpha B\cup C$ if $A\sim^{S_F}_\alpha B$.
\\{\bf Proof.} This result holds obviously for $A$ and $B$ being empty. For two nonempty fuzzy sets $A$ and $B$, we have  $S_F(A\cup C, B\cup C)=
\mathop{\inf}\limits_{x\in U}F((A\cup C)(x),(B\cup C)(x))=\mathop{\inf}\limits_{x\in U}F(A(x)\vee C(x),B(x)\vee C(x))=\mathop{\inf}\limits_{x\in U}M(I(A(x)\vee C(x),B(x)\vee C(x)),I(B(x)\vee C(x),A(x)\vee C(x)))$. Notice that $I(A(x)\vee C(x),B(x)\vee C(x))=I(A(x),B(x)\vee C(x))\wedge I(C(x),B(x)\vee C(x))=I(A(x),B(x)\vee C(x))\wedge 1=I(A(x),B(x)\vee C(x))\geq I(A(x),B(x))$. Similarly, we can obtain $I(B(x)\vee C(x),A(x)\vee C(x))=I(B(x),A(x)\vee C(x))\geq I(B(x),A(x))$. Since $M$ is non-decreasing in each variable,  $S_F(A\cup C, B\cup C)=\mathop{\inf}\limits_{x\in U}F(A(x)\vee C(x),B(x)\vee C(x))=\mathop{\inf}\limits_{x\in U}M(I(A(x)\vee C(x),B(x)\vee C(x)),I(B(x)\vee C(x),A(x)\vee C(x)))\geq\mathop{\inf}\limits_{x\in U}M(I(A(x),B(x),I(B(x),A(x)))=S_F(A, B)$ holds. Therefore, we have $S_F(A\cup C, B\cup C)\geq \alpha$. That is, $A\cup C\sim^{S_F}_\alpha B\cup C$.
\\{\bf Proposition 4.4} $A\cap C\sim^{S_F}_\alpha B\cap C$ if $A\sim^{S_F}_\alpha B$.
\\{\bf Proof.} It is easy to see that this statement holds when $A$ and $B$ are empty. For two nonempty fuzzy sets $A$ and $B$, we have $S_F(A\cap C, B\cap C)=
\mathop{\inf}\limits_{x\in U}F((A\cap C)(x),(B\cap C)(x))=\mathop{\inf}\limits_{x\in U}F(A(x)\wedge C(x),B(x)\wedge C(x))=\mathop{\inf}\limits_{x\in U}M(I(A(x)\wedge C(x),B(x)\wedge C(x)),I(B(x)\wedge C(x),A(x)\wedge C(x)))$. Notice that $I(A(x)\wedge C(x),B(x)\wedge C(x))=I(A(x)\wedge C(x),B(x))\wedge I(A(x)\wedge C(x),C(x))=I(A(x)\wedge C(x),B(x))\wedge 1=I(A(x)\wedge C(x),B(x))\geq I(A(x),B(x))$. We can similarly obtain $I(B(x)\wedge C(x),A(x)\wedge C(x))=I(B(x)\wedge C(x),A(x))\geq I(B(x),A(x))$. Since $M$ is non-decreasing in each variable, we have $S_F(A\cap C, B\cap C)=\mathop{\inf}\limits_{x\in U}F(A(x)\wedge C(x),B(x)\wedge C(x))=\mathop{\inf}\limits_{x\in U}M(I(A(x)\wedge$\vspace{1mm} $ C(x),B(x)\wedge C(x)),I(B(x)\wedge C(x),A(x)\wedge C(x)))\geq\mathop{\inf}\limits_{x\in U}M(I(A(x),B(x),I(B(x),A(x)))=$\vspace{1mm} $S_F(A, B)$. Thus, $S_F(A\cap C, B\cap C)\geq \alpha$ holds  if $S_F(A,B)\geq \alpha$. That is, $A\cap C\sim^{S_F}_\alpha B\cap C$.

Considering that Raha's SBAR method  is based on fuzzy relations, we need to investigate the retainment of approximate equality when composing a fuzzy relation with a fuzzy set.
\\{\bf Proposition 4.5} Let $R$ be a fuzzy relation from $U$ to $V$ and $R(A)(y)=\mathop{\sup}\limits_{x\in U} A(x)\wedge R(x, y)$.\vspace{1mm} We have $R(A) \sim^{S_F}_\alpha R(B)$ if the fuzzy implication $I$  is continuous in the first variable\vspace{1mm} and  $A\sim^{S_F}_\alpha B$.
\\{\bf Proof.} Obviously, this result holds when $A$ and $B$ are empty. For two nonempty
fuzzy sets $A$\vspace{1mm} and $B$, we have $S_F(R(A), R(B))=
\mathop{\inf}\limits_{y\in V}F(R(A),R(B))=\mathop{\inf}\limits_{y\in V}F\left(\mathop{\sup}\limits_{x\in U} A(x)\wedge R(x, y),\mathop{\sup}\limits_{x\in U}B(x)\right.$\vspace{1mm} $\wedge R(x, y))=\mathop{\inf}\limits_{y\in V}M(I\left(\mathop{\sup}\limits_{x\in U} A(x)\wedge R(x, y),\mathop{\sup}\limits_{x\in U}B(x)\wedge R(x, y)\right),I\left(\mathop{\sup}\limits_{x\in U}B(x)\wedge R(x, y),\mathop{\sup}\limits_{x\in U}A(x)\right.$\vspace{1mm} $\wedge R(x, y)))$. Since $I$ fulfills I2 and continuity for the first variable, we have $I\left(\mathop{\sup}\limits_{x\in U} A(x)\wedge R(x, y),\right.$\vspace{1mm} $\mathop{\sup}\limits_{x\in U}B(x)\wedge R(x, y))\geq \mathop{\sup}\limits_{x_1\in U}\mathop{\inf}\limits_{x_2\in U} I(A(x_2)$\vspace{1mm} $\wedge R(x_2, y),B(x_1)\wedge R(x_1, y))=\mathop{\sup}\limits_{x_1\in U}\mathop{\inf}\limits_{x_2\in U} (I(A(x_2)\wedge R(x_2, y),B(x_1))$\vspace{1mm} $\wedge I(A(x_2)\wedge R(x_2, y),R(x_1, y))$ $=\mathop{\sup}\limits_{x_1\in U}\mathop{\inf}\limits_{x_2\in U} (I(A(x_2)\wedge R(x_2, y),B(x_1))\wedge 1)\geq \mathop{\inf}\limits_{x\in U}I(A(x),B(x))$.\vspace{1mm} Similarly, $I\left(\mathop{\sup}\limits_{x\in U} A(x)\wedge R(x, y),\mathop{\sup}\limits_{x\in U}B(x)\wedge R(x, y)\right)\geq \mathop{\inf}\limits_{x\in U}I(B(x),A(x))$ holds. Since $M$ is non-decreasing in each variable, we obtain $S_F(R(A), R(B))=\mathop{\inf}\limits_{y\in V}F(R(A),$\vspace{1mm} $R(B))=\mathop{\inf}\limits_{y\in V}F\left(\mathop{\sup}\limits_{x\in U} A(x)\wedge R(x, y),\mathop{\sup}\limits_{x\in U}B(x)\wedge R(x, y)\right)\geq\mathop{\inf}\limits_{x\in U}M(I(A(x),B(x),I(B(x),A(x)))$\vspace{1.5mm} $=S_F(A, B)$. And then $S_F(R(A), R(B))\geq \alpha$ holds  if $S_F(A,B)\geq \alpha$.
\\{\bf Remark 6.} Indeed, the mapping $R(A)$ in Proposition 4.5 corresponds to the Compositional Rule of Inference (CRI) method introduced in \cite{Zadeh}. An arbitrary  t-norm can be applied in CRI method. Let us denote $\widetilde{R}(A)(y)=\mathop{\sup}\limits_{x\in U} T(A(x), R(x, y))$. Nevertheless, the following example demonstrates that $\widetilde{R}(A) \sim^{S_F}_\alpha \widetilde{R}(B)$\vspace{1mm} does not hold although $A\sim^{S_F}_\alpha B$.\vspace{1mm}
\\{\bf Example 4.6} Let $A=[0.7\ 0.8\ 0.4]$, $R=\left(
                                                     \begin{array}{ccc}
                                                       0.2&0.1&0.2\\
                                                       0.1&0.4&0.3\\
                                                        0.5&0.3&0.5\\
                                                     \end{array}
                                                   \right)$ and $T(x,y)=\left\{
                                                     \begin{array}{ll}
                                                       0&x+y<1\\
                                                       x\wedge y&x+y\geq 1\\
                                                     \end{array}
                                                   \right.$.\vspace{1mm}
 We have $\widetilde{R}(A)=[0\ 0.4\ 0.3]$. Further, suppose that
$A'=[0.9\ 0.6\ 0.6]$. Then, $\widetilde{R}(A')=[0.5\ 0.4\ 0.5]$.  Let $T_M$ and $I_L$ be used to describe $S_F$, we have $S_F(A,A')=0.8$. Take $\alpha=0.8$. However, $S_F(\widetilde{R}(A),\widetilde{R}(A'))=0.5<\alpha$.

As shown in Example 4.6, the approximation equality of $\widetilde{R}(A)$ cannot be retained. To analyze the approximation equality of $B'_1$ in Eq.(1) and $B_2'$ in Eq.(2), we impose the additional constraint that $F$ is generated by $M=T$ and $I=I_T$ in the rest of this work, that is, $F(x,y)=T(I_T(x,y),I_T(y,x))$. Consequently, the similarity measure $S_F$ admits the reformulation $S_F(A,B)=\mathop{\inf}\limits_{x\in U}T(I_T(A(x),B(x)),I_T(B(x),$ $A(x)))$.  We then have the following results.
\\{\bf Proposition 4.7} Let $T$ be a left continuous t-norm. If $A\sim^{S_F}_\alpha A'$ and $B\sim^{S_F}_\alpha B'$, then $T(A,B)\sim^{S_F}_{T(\alpha,\alpha)} T(A',B')$.
\\{\bf Proof.} It is sufficient to ensure that $S_F(T(A,B), T(A',B'))\geq T(S_F(A,A'),S_F(B,B'))$ holds.\vspace{1mm} As $I_T$ satisfies OP,  $S_F(A,A')=\mathop{\inf}\limits_{x\in U}F(A(x), A'(x))=\mathop{\inf}\limits_{x\in U}T(I_T(A(x), A'(x)), I_T(A'(x), $ $A(x)))=\mathop{\inf}\limits_{A(x)>A'(x)}I_T(A(x),A'(x))\wedge \mathop{\inf}\limits_{A(x)\leq A'(x) }I_T(A'(x),A(x))$ holds. Therefore, we can rewrite $T(S_F(A,A'),S_F(B,B'))=T\left(\mathop{\inf}\limits_{x\in U}F(A(x),A'(x)),\mathop{\inf}\limits_{x\in U}F(B(x),B'(x))\right)=T\left(\mathop{\inf}\limits_{x\in U}T(I_T(A(x),\right.$\vspace{1mm} $\left.A'(x)),I_T(A'(x), A(x))),\mathop{\inf}\limits_{x\in U}T(I_T(B(x), B'(x)),I_T(B'(x), B(x)))\right)=T\left(\mathop{\inf}\limits_{A(x)>A'(x)}I_T(A(x),\right.$\vspace{1mm} $\left.A'(x))\wedge \mathop{\inf}\limits_{A(x)\leq A'(x)}I_T(A'(x),A(x)),\ \mathop{\inf}\limits_{B(x)>B'(x)}I_T(B(x),B'(x))\wedge\mathop{\inf}\limits_{B(x)\leq B'(x)}I_T(B'(x),B(x))\right)=$\vspace{1.5mm} $\mathop{\inf}\limits_{A(x)>A'(x)}\mathop{\inf}\limits_{B(x)>B'(x)}T(I_T(A(x),A'(x)),I_T(B(x), B'(x)))\wedge \mathop{\inf}\limits_{A(x)>A'(x)}\mathop{\inf}\limits_{B(x)\leq B'(x)}T(I_T(A(x),A'(x)),$\vspace{1.5mm} $I_T(B'(x), B(x)))\wedge \mathop{\inf}\limits_{A(x)\leq A'(x)}\mathop{\inf}\limits_{B(x)> B'(x)}T(I_T(A'(x),A(x)),I_T(B(x), B'(x)))\wedge \mathop{\inf}\limits_{A(x)\leq A'(x)}\mathop{\inf}\limits_{B(x)\leq B'(x)}$\vspace{1.5mm} $T(I_T(A'(x),A(x)),I_T(B'(x), B(x)))$. On the other hand, We have $S_F(T(A,B), T(A',B'))=\mathop{\inf}\limits_{x\in U}F(T(A(x), B(x)), T(A(x), B(x)))=\mathop{\inf}\limits_{x\in U}T(I_T(T(A(x), B(x)), T(A'(x),B'(x))), I_T(T(A'(x),$\vspace{1mm}  $B'(x)), T(A(x), B(x))))=\mathop{\inf}\limits_{A(x)>A'(x)}\mathop{\inf}\limits_{B(x)>B'(x)}I_T(T(A(x), B(x)), T(A'(x),B'(x)))\wedge \mathop{\inf}\limits_{A(x)>A'(x)}$\vspace{1mm} $\mathop{\inf}\limits_{B(x)\leq B'(x)}T(I_T(T(A(x), B(x)), T(A'(x),B'(x))), I_T(T(A'(x),B'(x)), T(A(x), B(x))))\wedge \mathop{\inf}\limits_{A(x)\leq A'(x)}$\vspace{1mm} $\mathop{\inf}\limits_{B(x)> B'(x)}T(I_T(T(A(x), B(x)), T(A'(x),B'(x))), I_T(T(A'(x),B'(x)), T(A(x), B(x))))\wedge \mathop{\inf}\limits_{A(x)\leq A'(x)}$\vspace{1mm} $\mathop{\inf}\limits_{B(x)\leq B'(x)}I_T(T(A'(x),B'(x)), T(A(x), B(x)))$. Let us study the following four cases:\vspace{1mm}

i. If $A(x)\leq A'(x)$ and $B(x)\leq B'(x)$. We then have $T(A(x),B(x))\leq T(A'(x),B'(x))$. Thus, $T(I_T(T(A(x), B(x)), T(A'(x),B'(x))), I_T(T(A'(x),B'(x)), T(A(x), B(x))))=I_T(T(A'(x),$ $B'(x)), T(A(x), B(x)))$. The left continuity of $T$ implies that $T(A'(x),T(B'(x),T(I_T(A'(x),$ $A(x)),I_T(B'(x), B(x)))))=T(T(A'(x),I_T(A'(x),A(x))),T(B'(x),I_T(B'(x), B(x)))))\leq T(A(x),$ $B(x)))$ holds. Therefore, we have $T(I_T(A'(x), A(x)),I_T(B'(x), B(x)))))\leq I_T(T(A'(x),B'(x)),$ $T(A(x), B(x)))$ as $T$ is left continuous.

ii. If $A(x)\leq A'(x)$ and $B(x)>B'(x)$. In this case we cannot compare $T(A(x), B(x))$ with $T(A'(x), B'(x))$. Thus, we further consider the following two options.

ii-1. $T(A(x), B(x))\leq T(A'(x), B'(x))$. This means $T(I_T(T(A(x), B(x)), T(A'(x),B'(x))),$ $ I_T(T(A'(x),B'(x)), T(A(x), B(x))))= I_T(T(A'(x),B'(x)), T(A(x), B(x))))$. Then, the left continuity of $T$ implies that $T(A'(x),T(B'(x),T(I_T(A'(x),A(x)),I_T(B(x), B'(x)))))=T(T(A'(x),$ $I_T(A'(x),A(x))),T(B'(x),I_T(B(x), B'(x)))))\leq T(A(x),T(B'(x),I_T(B(x), B'(x))))$ holds. Since $I_T(B(x),B'(x))\leq I_T(B'(x),B(x))=1$, we have $T(B'(x),I_T(B(x),B'(x)))\leq B(x)$ by the left continuity of $T$. This implies that $T(A'(x),T(B'(x),T(I_T(A'(x),A(x)),I_T(B(x), B'(x)))))\leq T(A(x),T(B'(x),I_T(B(x), B'(x))))\leq T(A(x), B(x))$ holds. Thus, we have $T(I_T(A'(x),A(x)),$ $I_T(B(x), B'(x)))\leq I_T(T(A'(x),B'(x)), T(A(x), B(x)))$.

ii-2. $T(A(x), B(x))> T(A'(x), B'(x))$. Similar to ii-1, $T(I_T(A'(x),A(x)),I_T(B(x),$ $ B'(x)))\leq I_T(T(A(x), B(x)),T(A'(x),B'(x)))$ holds.

iii. If $A(x)>A'(x)$ and $B(x)\leq B'(x)$. This case can be investigated similarly to ii.

iv. If $A(x)> A'(x)$ and $B(x)> B'(x)$. Similar to i, we  obtain $T(I_T(A(x), A'(x)),I_T(B(x), $ $B'(x)))))\leq I_T(T(A(x),B(x)),T(A'(x), B'(x)))$.

In a word, we obtain $S_F(T(A,B), T(A',B'))=\mathop{\inf}\limits_{x\in U}T(I_T(T(A(x), B(x)), T(A'(x),B'(x))),$ $I_T(T(A'(x),B'(x)), T(A(x), B(x))))\geq T\left(\mathop{\inf}\limits_{x\in U}F(A(x),A'(x)),\mathop{\inf}\limits_{x\in U}F(B(x),B'(x))\right)=T(S_F(A,$ $A'),S_F(B,B'))$.
\\{\bf Proposition 4.8} Let $I_T$ be an R-implication generated by a left continuous t-norm $T$. If $A\sim^{S_F}_\alpha A'$ and $B\sim^{S_F}_\alpha B'$, then $I_T(A,B)\sim^{S_F}_{T(\alpha,\alpha)} I_T(A',B')$.
\\{\bf Proof.} We have $S_F(I_T(A,B), I_T(A',B'))=\mathop{\inf}\limits_{x\in U}F(I_T(A(x), B(x)), I_T(A(x), B(x)))=\mathop{\inf}\limits_{x\in U}T(I_T$\vspace{1mm} $(I_T(A(x),B(x)), I_T(A'(x),B'(x))), I_T(I_T(A'(x),B'(x)), I_T(A(x), B(x))))=\mathop{\inf}\limits_{A(x)\leq A'(x)}\mathop{\inf}\limits_{B(x)>B'(x)}$\vspace{1mm} $I_T(I_T(A(x),B(x)), I_T(A'(x),B'(x)))\wedge \mathop{\inf}\limits_{A(x)>A'(x)}\mathop{\inf}\limits_{B(x)\leq B'(x)}I_T(I_T(A'(x),B'(x)), I_T(A(x),B(x)))$\vspace{1mm} $\wedge \mathop{\inf}\limits_{A(x)> A'(x)}\mathop{\inf}\limits_{B(x)> B'(x)}T(I_T(I_T(A(x),B(x)), I_T(A'(x),B'(x))), I_T(I_T(A'(x),B'(x)), I_T(A(x),B($\vspace{1mm} $x))))\wedge \mathop{\inf}\limits_{A(x)\leq A'(x)}$\vspace{1mm} $\mathop{\inf}\limits_{B(x)\leq B'(x)}T(I_T(I_T(A(x),B(x)), I_T(A'(x),B'(x))), I_T(I_T(A'(x),B'(x)), I_T(A(x),$ $B(x))))$. Let us study the four cases as follows:

i. If $A(x)\leq A'(x)$ and $B(x)> B'(x)$. As $T$ is left continuous, we have $T(T(T(A'(x),I_T(A'(x),$ $A(x))), I_T(A(x), B(x))), I_T(B(x),B'(x)))\leq B'(x)$. Again, $T(I_T(A'(x),A(x))), I_T(B(x), B(x)))$ $\leq  I_T(I_T(I_T(A(x),B(x))), I_T(I_T(A'(x),B'(x))))$ holds.

ii. If $A(x)> A'(x)$ and $B(x)\leq B'(x)$. Similar to i, we can get  $T(I_T(A(x),A'(x))), I_T(B'(x),$ $B(x)))\leq  I_T(I_T(I_T(A'(x),B'(x))), I_T(I_T(A(x),B(x))))$.

iii. If $A(x)> A'(x)$ and $B(x)> B'(x)$. We further study the following two options.

iii-1. $I_T(A(x), B(x))\leq I_T(A'(x), B'(x))$. This case means $T(I_T(I_T(A(x), B(x)), I_T(A'(x),$ $ B'(x))),I_T(I_T(A'(x),B'(x)), I_T(A(x), B(x))))= I_T(I_T(A'(x),B'(x)), I_T(A(x), B(x))))$. Then, the left continuity of $T$ implies that $T(T(T(A(x),I_T(A(x),A'(x))), I_T(A'(x), B'(x))), I_T(B(x),$ $B'(x)))\leq T(B'(x), I_T(B(x),B'(x)))$ holds. Since $I_T(B(x),B'(x))\leq I_T(B'(x),B(x))=1$, we have $T(B'(x),I_T(B(x),B'(x)))\leq B(x)$ by the left continuity of $T$. This implies that $T(T(I_T(A(x),A'(x)), I_T(A'(x), B'(x))), I_T(B(x),$ $B'(x)))\leq I_T(A(x),B(x)))$ holds. Thus, we have $T(I_T(A'(x),A(x)),$ $I_T(B(x), B'(x)))\leq I_T(I_T(A'(x),B'(x)), I_T(A(x), B(x)))$.

iii-2. $I_T(A(x), B(x))> I_T(A'(x), B'(x))$. Similar to iii-1, $T(I_T(A(x),A'(x)),I_T(B(x),$ $B'(x)))\leq I_T(T(A(x), B(x)),T(A'(x),B'(x)))$ holds.

iv. If $A(x)\leq A'(x)$ and $B(x)\leq B'(x)$. Similar to iii, we  obtain $T(I_T(A(x), A'(x)),I_T(B(x), $ $B'(x)))))\leq I_T(T(A(x),B(x)),T(A'(x), B'(x)))$.

In a word, this means that $S_F(I_T(A,B),I_T(A',B'))\geq T(S_F(A,A'),S_F(B,B'))$ holds.
\\{\bf Proposition 4.9} $|S_F(A,B)-S_F(A',B')|\leq 2-S_F(A,A')-S_F(B,B')$.
\\{\bf Proof.} It is enough to ensure that $F(x,y)+F(y,z)\leq F(x,z)+1$ holds for any $x,y,z\in[0,1]$. Without loss of generality, let $x\leq y$. We further consider the following three options:

i. If $x\leq y\leq z$. This implies that $F(x,y)=I(y,x)$, $F(y,z)=I(z,y)$ and $F(x,z)=I(z,x)$. Therefore, $F(x,y)+F(y,z)=I(y,x)+I(z,y)\leq I(z,x)+1= F(x,z)+1$.

ii. If $x\leq z\leq y$. We then have $F(x,y)=I(y,x)$, $F(y,z)=I(y,z)$ and $F(x,z)=I(z,x)$. This means that $F(x,y)+F(y,z)=I(y,x)+I(y,z)\leq I(z,x)+1= F(x,z)+1$ holds.

iii. If $z\leq x\leq y$. This case means that $F(x,y)=I(y,x)$, $F(y,z)=I(y,z)$ and $F(x,z)=I(x,z)$. We therefore obtain $F(x,y)+F(y,z)=I(y,x)+I(y,z)\leq I(x,z)+1= F(x,z)+1$.\vspace{1mm}
\\{\bf Proposition 4.10} Let $f$ and $g$ be two mappings from $U$ to [0,1]. We have $F\left(\mathop{\inf}\limits_{x\in U}f(x), \mathop{\inf}\limits_{x\in U}g(x)\right)$\vspace{1mm} $\geq\mathop{\inf}\limits_{x\in U}F(f(x),g(x))$ and $F\left(\mathop{\sup}\limits_{x\in U}f(x), \mathop{\sup}\limits_{x\in U}g(x)\right)\geq\mathop{\inf}\limits_{x\in U}F(f(x),g(x))$.\vspace{1mm}
\\{\bf Proof.} We only prove the inequality $F\left(\mathop{\inf}\limits_{x\in U}f(x), \mathop{\inf}\limits_{x\in U}f(x)\right)\geq\mathop{\inf}\limits_{x\in U}F(f(x),g(x))$. The other\vspace{1mm} can be proved similarly. We have $\mathop{\inf}\limits_{x\in U}F(f(x),g(x))=\mathop{\inf}\limits_{f(x)>g(x)}I_T(f(x),g(x))\wedge \mathop{\inf}\limits_{f(x)\leq g(x)}I_T(g(x),$ $f(x))$. Let us consider the two cases in the following:

i. If $f(x)>g(x)$. This implies that $\mathop{\inf}\limits_{x\in U}f(x)\geq \mathop{\inf}\limits_{x\in U}g(x)$. And then $F\left(\mathop{\inf}\limits_{x\in U}f(x), \mathop{\inf}\limits_{x\in U}f(x)\right)=$\vspace{1mm} $I_T\left(\mathop{\inf}\limits_{x\in U}f(x), \mathop{\inf}\limits_{x\in U}g(x)\right)=\mathop{\inf}\limits_{x\in U}I_T\left(\mathop{\inf}\limits_{x\in U}f(x),g(x)\right)\geq \mathop{\inf}\limits_{f(x)>g(x)}I_T(f(x),g(x))$. This means\vspace{1mm}  that $F\left(\mathop{\inf}\limits_{x\in U}f(x), \mathop{\inf}\limits_{x\in U}g(x)\right)\geq\mathop{\inf}\limits_{x\in U}F(f(x),g(x))$ holds.

ii. If $f(x)\leq g(x)$. This implies that $\mathop{\inf}\limits_{x\in U}f(x)\leq \mathop{\inf}\limits_{x\in U}g(x)$. And then $F\left(\mathop{\inf}\limits_{x\in U}f(x), \mathop{\inf}\limits_{x\in U}f(x)\right)$\vspace{1mm} $=I_T\left(\mathop{\inf}\limits_{x\in U}g(x), \mathop{\inf}\limits_{x\in U}f(x)\right)=\mathop{\inf}\limits_{x\in U}I_T\left(\mathop{\inf}\limits_{x\in U}g(x),f(x)\right)\geq \mathop{\inf}\limits_{f(x)>g(x)}I_T(f(x),g(x))$. This means\vspace{1mm}  that $F\left(\mathop{\inf}\limits_{x\in U}f(x), \mathop{\inf}\limits_{x\in U}g(x)\right)\geq\mathop{\inf}\limits_{x\in U}F(g(x),f(x))$ holds.

In a word, we obtain $F\left(\mathop{\inf}\limits_{x\in U}f(x), \mathop{\inf}\limits_{x\in U}f(x)\right)\geq\mathop{\inf}\limits_{x\in U}F(f(x),g(x))$.\vspace{1mm}

Finally, we study the approximation equality of Raha's SBAR method with the preceding arguments. Let $B''_1(y) =\mathop{\sup}\limits_{x\in U} I_T(S_F(A'',A_1),$\vspace{1mm} $ T(A_1(x),B_1(y)))$. According to Proposition 4.9, we have $|S_F(A'',A_1)-S_F(A',A)|\leq \varepsilon\ (=2-\alpha_1-\alpha_2)$ when $A'\sim^{S_F}_{\alpha_1} A''$ and $A\sim^{S_F}_{\alpha_2} A_1$. Setting  $S_F(A',A)=\alpha$, the inequality  $\alpha-\varepsilon\leq S_F(A'',A_1)\leq \alpha+\varepsilon$ holds. Further, we have the following statement.
 \\{\bf Theorem 4.11} If $A'\sim^{S_F}_{\alpha_1} A''$, $A\sim^{S_F}_{\alpha_2} A_1$ and $B\sim^{S_F}_{\alpha_3} B_1$ in Eq.(1), then $B_1' \sim^{S_F}_{T(\beta,\alpha_2,\alpha_3)} B''_1$ with $\beta=\mathop{\inf}\limits_{x\in U}F(S_F(A'(x),A(x)),S_F(A''(x),A_1(x)))$.
\\{\bf Proof.} According to Propositions 4.7, 4.8 and 4.10, we have $S_F(B_1',B_1'')=\mathop{\inf}\limits_{y\in V}F(B_1'(y),B_1''(y))$\vspace{1mm} $=\mathop{\inf}\limits_{y\in V}F\left(\mathop{\sup}\limits_{x\in U} I_T(S_F(A',A),T(A(x),B(y))),\mathop{\sup}\limits_{x\in U} I_T(S_F(A'',A_1),T(A_1(x),B_1(y)))\right)$\vspace{1mm} $\geq \mathop{\inf}\limits_{y\in V}\mathop{\inf}\limits_{x\in U}$ $F( I_T(S_F(A',A),T(A(x),B(y))),I_T(S_F(A'',A_1),T(A_1(x),B_1(y))))\geq \mathop{\inf}\limits_{y\in V}$\vspace{1mm} $\mathop{\inf}\limits_{x\in U} T(F(S_F(A',A),$ $S_F(A'',A_1)),F(T(A(x),B(y)),T(A_1(x),B_1(y))))\geq \mathop{\inf}\limits_{y\in V}\mathop{\inf}\limits_{x\in U}T(F(S_F(A',A),S_F(A'',A_1)), T(F$ $(A(x),A_1(x)),T(B(y),B_1(y))))\geq T(\beta,\alpha_2, \alpha_2)$, where $\beta=\mathop{\inf}\limits_{x\in U}F(S_F(A'(x),A(x)),S_F(A''(x),$ $A_1(x)))$.\vspace{1mm}

Theorem 4.11 shows that the  approximation equality of $B_1'$ in Eq.(1) can be expressed by that of $A'$, $A$ and $B$. Similarly, let $B''_2(y) =\mathop{\inf}\limits_{x\in U} I_T(S(A'',A_1), I_T(A_1(x), B_1(y))$. We can obtain the  approximation equality of $B_2'$ in Eq.(2) as shown in the following.
\\{\bf Theorem 4.12} If $A'\sim^{S_F}_{\alpha_1} A''$, $A\sim^{S_F}_{\alpha_2} A_1$ and $B\sim^{S_F}_{\alpha_3} B_1$ in Eq.(2), $B_2' \sim^{S_F}_{T(\beta,\alpha_2,\alpha_3)} B_2''$ with $\beta=\mathop{\inf}\limits_{x\in U}F(S_F(A'(x),A(x)),S_F(A''(x),A_1(x)))$.\vspace{1mm}
\\{\bf Proof.} This can be proved similarly to Theorem 4.11.
\\{\bf Remark 7.} i. Obviously, $\beta$ can be expressed by $\alpha_1$ and $\alpha_2$;

ii. With Theorems 4.11 and 4.12, we find the fact that the closer $A'$ is to $A''$, the closer $B'$ is to $B$ when the antecedent and consequent of if-then rule are specified. Especially, $B'=B$ if $\alpha_1=\alpha_2=\alpha_3=1$, which aligns with human intuition.
\section{Hierarchical Raha's SBAR method with REF}
\qquad As discussed in Section 4, REF can be used to describe the similarity of two fuzzy sets. This enables the direct application of REF to SBAR method. To promote the computational efficiency of Raha's SBAR method, this section will construct two hierarchical Raha's SBAR methods with REF mentioned in Section 4. The computational and space complexity of these hierarchical Raha's SBAR methods are subsequently analyzed.

For methodological demonstration, we only consider the fuzzy system containing a single if-then rule, formally expressed as

IF $x_1$ is $A_{1}$ AND $x_2$ is $A_{2}$ THEN  $y$ is  $B$.

Let $A'=A'_1\times A'_2$ be a fuzzy input on the universe $U=U_1\times U_2$. If the similarity  between $A'$ and $A=A_1\times A_2$ is measured by $S_F(A',A)=T(S_F(A_1',A_1),S_F(A_2',A_2))$ (see Ref.\cite{Turksen}), then we have the following fact.
\\{\bf Theorem 5.1} If the logical connective``AND" in the if-then rule is interpreted by the t-norm $T$, then the output $B'_2$ in Eq.(2) is $B'_2(y)=\inf\limits_{x_1\in U_1}I_T(S_F(A'_1,A_1),I_T(A_{1}(x_1),B'_{21}(y)))$ with\vspace{1mm} $B'_{21}(y)=\inf\limits_{x_2\in U_2}I_T(S_F(A'_{2},A_2),I_T(A_{2}(x_2),B(y))))$.\vspace{1mm}
\\{\bf Proof.} Since $I_T$ satisfies EP and $I_T(T(x,y),z)=I_T(x,I_T(y,z))$ holds for the t-norm $T$ according to Ref.\cite{Baczynski}, we have  $I_T(T(S_F(A_1',A_1),S_F(A_2',A_2)),I_T(T(A_{1}(x_1),A_{2}(x_2)),B(y)))=I_T(S_F(A_1',A_1),$ $I_T(S_F(A_2',A_2),I_T(A_{1}(x_1),I_T(A_{2}(x_2),B(y)))))=I_T(S_F(A_1',A_1),I_T(A_{1}(x_1),I_T$ $(S_F(A_2',A_2),I_T(A_{2}(x_2),B(y)))))$. This means that
$B'_2$ in Eq.(2) can be rewritten as
$B'_2(y)=\inf\limits_{(x_1, x_2)\in U_1\times U_2}I_T$\vspace{1mm} $(T(S_F(A'_{1},A_1),S_F(A'_2,A_2)),I_T(T(A_{1}(x_1),A_{2j}(x_2)),B_j(y)))=\inf\limits_{(x_1, x_2)\in U_1\times U_2}I_T$\vspace{1mm} $(S_F(A_1',A_1),I_T(A_{1}(x_1),I_T(S_F(A_2',A_2),I_T(A_{2}(x_2),B(y)))))=\inf\limits_{x_1\in U_1}I_T(S_F(A_1',A_1),I_T(A_{1}(x_1),$\vspace{1mm} $\inf\limits_{x_2\in U_2}I_T(S_F(A_2',A_2),I_T(A_{2}(x_2),B(y)))))$. Let  $B'_{21}(y)=\inf\limits_{x_2\in U_2}I_T(S_F(A'_{2},A_2),I_T(A_{2}(x_2),B(y))))$.\vspace{1mm} Therefore, we obtain $B'_2(y)=\inf\limits_{x_1\in U_1}I_T(S_F(A'_1,A_1),I_T(A_{1}(x_1),B'_{21}(y)))$.\vspace{1mm}

Theorem 5.1 shows that the output of fuzzy system  can be derived through independent computations of $B'_{21}$ and $B'_2$. This stepwise calculation process characterizes the hierarchical Raha's SBAR method, which serves as the foundation for constructing Algorithm 1.

When handling $n$-dimensional fuzzy inputs in a fuzzy system, it is preferable to treat the composite input $A'=A'_1\times A'_2\times\cdots\times A'_n$ as an entirety\cite{Wang}. Here, $A'$ denotes the fuzzifier of an $n$-dimensional input $\mathbf{x}\in U$. The membership grade of $A'$ is  typically computed as $A'(\mathbf{x})=T(A'_1(x_1),A'_2(x_2),\cdots, A'_n(x_n))$, where  $T$ is a t-norm interpreting the logical connective ``AND" in the if-then rules\cite{Wang}. The hierarchical SBAR method requires the condition $S_F(A',A)=T(S_F(A_1',A_1),S_F(A_2',A_2))$ to hold, analogous to Theorem 5.1. That is, we need to seek a t-norm $T$ and a fuzzy implication $I$ such that $S_F(T(A_1',A_2'),T(A_1,A_2))=T(S_F(A_1',A_1),S_F(A_2',A_2))$ holds. To achieve this objective, we seek the condition under which the following equation is satisfied
\begin{equation}
I(T(x,y),T(x',y'))=T(I(x,x'),I(y,y')),\ \forall\ x>x',y> y'.
 \end{equation}
 Thus, we suppose that $T$ is left continuous and $I=I_T$ again. The following counterexample demonstrates that left-continuity of $T$ alone is insufficient to guarantee the validity of Eq.(3).\vspace{1mm}
 \\{\bf Example 5.2} Let $T(x,y)=\left\{\begin{array}{ll}
                                               0& x+y\leq 1 \\
                                               x\wedge y&\textmd{otherwise}
                                             \end{array}\right.$. Then, the R-implication $I_T$ generated by\vspace{1mm} $T$ is $I_T(x,y)=\left\{\begin{array}{ll}
                                               1& x\leq y \\
                                               (1-x)\wedge y&\textmd{otherwise}
                                             \end{array}\right.$\cite{Klement}.\vspace{1mm}
 Taking $x=0.8,x'=0.5,y=0.4$ and $y'=0.3$, we have $I_T(T(0.8,0.6),T(0.4,0.3))=0.5\neq 0.6=T(I_T(0.8,0.4),I_T(0.6,0.3))$.

We therefore turn to continuous t-norms to determine the conditions under which Eq.(3) holds.
 \\{\bf Proposition 5.3} $T_M$ and $I_G$ satisfy Eq.(3).
\\{\bf Proof.} This can be directly verified.
\\{\bf Proposition 5.4} Let $T$ be a strict t-norm. $T$ and $I_T$ satisfies Eq.(3).
\\{\bf Proof.} Since $T$ is a strict t-norm, we can find an increasing bijection $t$ on [0,1] such that\vspace{1mm} $T(x,y)=t^{-1}(t(x)t(y))$\cite{Klement}. In this case,  $I_T(x,y)=\left\{\begin{array}{ll}
                                               1& x\leq y \\
                                               t^{-1}\left(\frac{t(x)}{t(y)}\right)&\textmd{otherwise}
                                             \end{array}\right.$\cite{Baczynski}. Thus,\vspace{1mm} for any $x>x'$ and $y>y'$, we have
$T(I_T(x,x'),I_T(y,y'))=T\left(t^{-1}\left(\frac{t(x')}{t(x)}\right),  t^{-1}\left(\frac{t(y')}{t(y)}\right)\right)=$\vspace{1mm} $t^{-1}\left(t\left(t^{-1}\left(\frac{t(x')}{t(x)}\right)\right)t\left(t^{-1}\left(\frac{t(y')}{t(y)}\right)\right)\right)=t^{-1}\left(\frac{t(x')}{t(x)}\frac{t(y')}{t(y)}\right)$. On the other hand,
$I_{T}(T(x,y),T(x',y'))$\vspace{1.5mm} $=t^{-1}\left(\frac{t(T(x',y'))}{t(T(x,y))}\right)=t^{-1}\left(\frac{t(t^{-1}(t(x')t(y')))}{t(t^{-1}(t(x)t(y)))}\right)=
t^{-1}\left(\frac{t(x')t(y'))}{t(x)t(y))}\right)$. This means that $T(I_T(x,x'),I_T(y,$\vspace{1mm} $y'))=I_{T}(T(x,y),T(x',y'))$ holds for any $x\geq x'$ and $y\geq y'$.
\\{\bf Proposition 5.5} Let $T$ be a nilpotent t-norm. If neither $x,y$ nor $x',y'$ are zero divisors of $T$, then $T$ and $I_T$ satisfies Eq.(3).
\\{\bf Proof.} As $T$ is nilpotent, there exists an increasing bijection $t$ on [0,1] such that $T(x,y)=$\vspace{1mm} $t^{-1}((t(x)+t(y)-1)\vee0)$\cite{Klement}. In this case, $I_T$ can be expressed as $I_T(x,y)=t^{-1}((1-t(x)+t(y))\wedge 1)$\cite{Baczynski}.\vspace{1mm} For any $x\geq x'$ and $y\geq y'$, we then obtain $T(I_T(x,x'),I_T(y,y'))=T(t^{-1}((1-t(x)+t(x'))\wedge 1),t^{-1}((1-t(y)+t(y'))\wedge 1))= t^{-1}((t(t^{-1}((1-t(x)+t(x'))\wedge 1))+t(t^{-1}((1-t(y)+t(y'))\wedge 1))-1)\vee 0)=t^{-1}((1-t(x)+t(x'))\wedge 1+(1-t(y)+t(y'))\wedge 1)-1)\vee 0)=t^{-1}((1-t(x)+t(x')+1-t(y)+t(y')-1)\vee 0)=t^{-1}(1-t(x)-t(y)+t(x')+t(y'))$. On the other hand,
$I_{T}(T(x,y),T(x',y'))=t^{-1}(1-t(T(x,y))+t(T(x',y')))=t^{-1}(1-t(t^{-1}((t(x)+t(y)-1)\vee0))+t(t^{-1}((t(x')+t(y')-1)\vee0)))=
t^{-1}(1-(t(x)+t(y)-1)\vee0)+(t(x')+t(y')-1)\vee0))$. Since neither $x,y$ nor $x',y'$ are zero divisors of $T$, we have $t(x)+t(y)>1$ and $t(x')+t(y')-1$. This means that $I_{T}(T(x,y),T(x',y'))=t^{-1}(1-t(x)-t(y)+t(x')+t(y'))=T(I_T(x,x'),I_T(y,y'))$ holds.
\\{\bf Remark 8.} For a  nilpotent t-norm $T$, both $x,y$ and $x',y'$ are not zero divisors of $T$ when Eq.(3) is satisfied for all $x>x'$ and $y>y'$. Otherwise, let $x,y$ be zero divisors of $T$. We then have $I_{T}(T(x,y),T(x',y'))=1$. On the other hand, $T(I_T(x,x'),I_T(y,y'))=t^{-1}(1-t(x)-t(y)+t(x')+t(y'))<1=I_{T}(T(x,y),T(x',y'))$. This contradicts with Eq.(3). Similarly, we can obtain a contradiction when $x',y'$ are zero divisors of $T$. This means that Eq.(3) holds if and only if neither $x,y$ nor $x',y'$ are zero divisors of $T$.
\\{\bf Proposition 5.6} Let $T=T_M$. For any two fuzzy inputs $A_1'$ and $A_2'$, $S_F(T(A_1,A_2),T(A_1',$ $A_2'))=T(S_F(A_1,A_1'),S_F(A_2,A_2'))$ holds  given the if-then rule with antecedents  $A_1$ and $A_2$.
\\{\bf Proof.} i. Let $T=T_M$.  As discussion in Section 4,  $S_F(T(A_1,A_2),T(A_1',A_2'))$ can be expressed as $S_F(A_1\wedge A_2,A_1'\wedge A_2')=\mathop{\inf}\limits_{A_1(x)>A_1'(x)}\mathop{\inf}\limits_{A_2(x)>A_2'(x)}(A_1'(x)\wedge A_2'(x))\wedge \mathop{\inf}\limits_{A_1(x)>A_1'(x)}\mathop{\inf}\limits_{A_2(x)\leq A_2'(x)}I_G(A(x)$\vspace{1mm} $\wedge A_2(x)), A_1'(x)\wedge A_2'(x))\wedge I_G(A_1'(x)\wedge A_2'(x), A_1(x)\wedge A_2(x))\wedge \mathop{\inf}\limits_{A_1(x)\leq A_1'(x)}\mathop{\inf}\limits_{A_2(x)> A_2'(x)}I_G(A_1(x)\wedge A_2(x), A_1'(x)$\vspace{1mm} $\wedge A_2'(x))\wedge I_G(A_1'(x)\wedge A_2'(x), A_!(x)\wedge A_2(x)))\wedge \mathop{\inf}\limits_{A_1(x)\leq A_1'(x)}\mathop{\inf}\limits_{A_2(x)\leq A_2'(x)}$\vspace{1mm} $(A_1(x)\wedge A_2(x))$. Further, we can obtain  $I_G(A_1(x)\wedge A_2(x)), A_1'(x)\wedge A_2'(x))\wedge I_G(A_1'(x)\wedge A_2'(x), A_1(x)\wedge A_2(x))=A_1'(x)\wedge A_2'(x)\wedge A_1(x)\wedge A_2(x)=A_1'(x)\wedge A_2(x)$. Similarly, $I_G(A_1(x)\wedge A_2(x), A_1'(x)\wedge A_2'(x))\wedge I_G(A_1'(x)\wedge B_2'(x), A_1(x)\wedge A_2(x))=A_1(x)\wedge A_2'(x)$ holds. Thus, we have $S_F(A_1\wedge A_2,A_1'\wedge A_2')=\mathop{\inf}\limits_{A_1(x)>A_1'(x)}\mathop{\inf}\limits_{A_2(x)>A_2'(x)}(A_1'(x)\wedge A_2'(x))\wedge$\vspace{1mm} $\mathop{\inf}\limits_{A_1(x)>A_1'(x)}\mathop{\inf}\limits_{A_2(x)\leq A_2'(x)}(A_1'(x)\wedge A_2(x))\wedge \mathop{\inf}\limits_{A_1(x)\leq A_1'(x)}\mathop{\inf}\limits_{A_2(x)> A_2'(x)}A_1(x)\wedge A_2'(x))\wedge \mathop{\inf}\limits_{A_1(x)\leq A_1'(x)}\mathop{\inf}\limits_{A_2(x)\leq A_2'(x)}$\vspace{1mm} $(A_1(x)\wedge A_2(x))=S_F(A_1,A_1')\wedge S_F(A_2,A_2')$.
\\{\bf Proposition 5.7} Let $T$ be a continuous Archimedean t-norm.  For any two fuzzy inputs $A_1'$ and $A_2'$, if either of the following conditions is satisfied

i. $A_1'\subseteq A_1$ and $A_2'\subseteq A_2$,

ii. $A_1\subseteq A_1'$ and $A_2\subseteq A_2'$,
\\ then $S_F(T(A_1,A_2),T(A_1',$ $A_2'))=T(S_F(A_1,A_1'),S_F(A_2,A_2'))$ holds given the if-then rule with antecedents.
\\{\bf Proof.} Let $T$ be a continuous Archimedean t-norm. Thus, $T$ is nilpotent or strict. Without loss of generality, we suppose that $T$ is a nilpotent
t-norm. It is sufficient to ensure that $S_F(T_L(A_1,$ $ A_2)),T_L(A_1',A_2'))=T_L(S_F(A_1,A_1'),S_F(A_2,A_2'))$ holds for $T_L$ and $I_L$.
Let us further investigate the following two cases:

i. If $A_1'\subseteq A_1$ and $A_2'\subseteq A_2$. $S_F(T_L(A_1,A_2)),$ $T_L(A_1',A_2'))$ can be expressed as $S_F(T_L(A_1,$ $ A_2),T_L(A_1',A_2'))=\mathop{\inf}\limits_{A_1(x)>A_1'(x)}\mathop{\inf}\limits_{A_2(x)>A_2'(x)}I_L(T_L(A_1(x), A_2(x)), T_L(A_1'(x),A_2'(x)))$. Since $A_1'$\vspace{1mm} and $A_2'$ are fuzzy inputs, $T_L(A_1'(x),A_2'(x))\neq 0$ for some $x\in U$. According to Proposition 5.3, we have $S_F(T_L(A_1, A_2), T_L(A_1',A_2'))=\mathop{\inf}\limits_{A_1(x)>A_1'(x)}\mathop{\inf}\limits_{A_2(x)>A_2'(x)}T_L(I_L(A_1(x),A_1'(x)),I_L(A_2(x), $\vspace{1mm} $A_2'(x)))= T_L(S_F(A_1,A_1'),S_F(A_2,A_2'))$.

ii. If $A_1\subseteq A_1'$ and $A_2\subseteq A_2'$. This case implies that $S_F(T_L(A_1, A_2),$ $T_L(A_1',A_2'))=\mathop{\inf}\limits_{A_1(x)\leq A_1'(x)}\mathop{\inf}\limits_{A_2(x)\leq A_2'(x)}I_L(T_L$\vspace{1mm} $(A_1'(x),A_2'(x)), T_L(A_1(x), A_2(x)))$. As $A_1$ and $A_2$ are the antecedents of the if-then rule, $T_L(A_1(x),A_2(x))\neq 0$ for some $x\in U$. By Proposition 5.5, we have $S_F(T_L(A_1, A_2),T_L(A_1',A_2'))= \mathop{\inf}\limits_{A_1(x)\leq A_1'(x)}\mathop{\inf}\limits_{A_2(x)\leq A_2'(x)}T_L(I_L(A_1'(x),A_1(x)),I_L(A_2'(x), A_2(x)))=$\vspace{1mm} $T_L(S_F(A_1,A_1'),S_F(A_2,A_2'))$.
 \\{\bf Remark 9.} The following counterexample shows that  the conditions i or ii of Proposition 5.7 is necessary.
\\{\bf Example 5.8} Let  $A_1=[0.9,0.6,0.7]$, $A_2=[0.4,0.6,0.5,0.3]$,  $A_1'= [1,0.5,0.8]$ and  $A_2'=[0.5,0.3,0.4,0.2]$. We have $S_F(T_L(A_1,A_2)),T_L(A_1',A_2'))=0.6$. However, $T_L(S_F(A_1,A_1'),S_F(A_2,$ $A_2'))=0.4$.

Indeed, it is not easy to ensure that the equation $S_F(T(A_1, A_2)),T(A_1',A_2'))=T(S_F(A_1,A_1'),$ $S_F(A_2,A_2'))$ holds. However, considering that $S_F(T(A_1, A_2)),T(A_1',A_2'))$ serves as the triggering  an if-then being fired, we can relax the original constraint to $S_F(T(A_1, A_2)),T(A_1',A_2'))\geq T(S_F(A_1,A_1'),S_F(A_2,A_2'))$. Similar to discussion in Section 4, $S_F(T(A_1, A_2)),T(A_1',A_2'))\geq T(S_F(A_1,A_1'),S_F(A_2,A_2'))$ holds for a continuous Archimedean t-norm. In this case, we can obtain the following result.
\\{\bf Theorem 5.9} Let the logical connective ``AND" in the if-then rule be interpreted by either the $T_M$ or a continuous Archimedean t-norm $T$. If the similarity  between $A'=A_1'\times A_2'$ and $A=A_1\times A_2$ is measured by $S_F$ which is generated by the same $T$, then the output $B'_2$ in Eq.(2) is $B'_2(y)=\inf\limits_{x_1\in U_1}I_T(S_F(A'_1,A_1),I_T(A_{1}(x_1),$\vspace{1mm} $B'_{21}(y)))$ with $B'_{21}(y)=\inf\limits_{x_2\in U_2}I_T(S_F(A'_{2},A_2),I_T(A_{2}(x_2),B(y))))$.\vspace{1mm}
 \\{\bf Proof.} This proof is similar to that of Theorem 5.1.

It is not difficult to see that this hierarchical Raha's SBAR method can be extended to compute the output of $n$-input-1-output fuzzy system. Let $U=U_1\times U_2\times\cdots\times U_n$. We further assume $|U_i|=u_i$ and $|V|=m$. Therefore,  $|U|=\mathop{\prod}\limits_{i=1}^n u_i$. And then an algorithm of this hierarchical Raha's SBAR method can be consequently provided as follows.
\\{\bf Algorithm 1} Hierarchical Raha's SBAR method with Eq.(3).

\textbf{Input:} $A_k'$, $A_k(k=1,2\cdots,n)$, $B$

\textbf{Output:} $B'$
1. \textbf{for} $i=1$ to $n$ and $j=2$ to $n$ \textbf{do}
\\ 2.  \quad\ \ $s_{i}=S_F(A_i',A_i)$
\\ 3. \quad\ \ $B'_{j}(y)=\inf\limits_{x_{i}\in U_{i}}I_T(s_{i},I_T(A_{i}(x_{i}),B(y)))$\vspace{1mm}
\\ 4. \quad\ \ $B\leftarrow B'_{j}$
\\ 5. \textbf{end for}
\\ 6. \textbf{return} $B'=B$

The computational and space complexity of Algorithm 1 can be investigated respectively in the following.
\begin{itemize}
  \item Computational complexity.  Obviously, the computational complexity of Algorithm 1 is $\mathcal{O}(m+u_i)$.
  \item Space complexity. In Algorithm 1, we only need to store $A_{i}(i=1,2,\cdots,n)$ and $B$. Thus, the space complexity of Algorithm 1 is $\mathcal{O}(m+u)$, where $u=\mathop{\max}\limits_{i=1}^nu_i$.
\end{itemize}
{\bf Remark 10.} i. It is easy to find that  the computational complexity of hierarchical Raha's SBAR method is a polynomial function of $n$. This polynomial-time characteristic ensures that the operation count in Algorithm 1 grows at a controlled rate, rather than exponentially, as the dimensionality of the if-then rules increases. Consequently, our proposed  hierarchical Raha's SBAR method can restrain the explosion of if-then rules.

ii. However, in Raha's SBAR method,  the computation of $S_F(A_1'\times\cdots\times A_n',A_1\times\cdots\times A_n)$ requires $\mathcal{O}(\mathop{\prod}\limits_{i=1}^n u_i+n)$ time. Similarly,  computing $ \inf\limits_{\mathbf{x}\in U}I_T(A(\mathbf{x}),B(y))$ has a time complexity of  $\mathcal{O}(\mathop{\prod}\limits_{i=1}^n u_i+m)$ time\cite{Cornelis}. The computational complexity of Raha's SBAR method s therefore given by $\mathcal{O}(\mathop{\prod}\limits_{i=1}^mu_i+n+m)$,which exhibits exponential dependence on the dimension $n$. Moreover, a $u$-dimensional matrix with $\mathop{\prod}\limits_{i=1}^mu_i$ entries needs to be stored for $A=A_1\times A_2\times\cdots\times A_n$\vspace{1mm} in Raha's SBAR
method. This implies the storage complexity of Raha's SBAR
method is\vspace{1mm} $\mathcal{O}(\mathop{\prod}\limits_{i=1}^mu_i+m)$\cite{Demirli}. Obviously, it is an\vspace{1mm} exponential function of $n$, too. Comparing Raha's SBAR
method with the hierarchical Raha's SBAR method, their computation and space complexity are listed together in Table 2.
 $$\mbox{\bf{\small Table\ 2 \ Comparison of computational and space complexity}}\vspace{-1mm}$$
\begin{center}
 \tabcolsep 0.05in
\begin{tabular}{ccc}
 \toprule[1pt]
 &Computational complexity&Space complexity\vspace{1mm}\\
\midrule[0.75pt]
Raha's SBAR method &$\mathcal{O}\left(\mathop{\prod}\limits_{i=1}^mu_i+n+m\right)$&$\mathcal{O}\left(\mathop{\prod}\limits_{i=1}^mu_i+m\right)$\vspace{1.5mm}\\
Hierarchical Raha's SBAR method&$\mathcal{O}(m+u_i)$&$\mathcal{O}(m+u)$\\
   \bottomrule[1pt]
\end{tabular}
\end{center}\vspace{1mm}
In a word, the hierarchical Raha's SBAR method has clear advantage in calculation and storage.

Finally, let us construct another hierarchical Raha's SBAR method with Eq.(1). Building on the methodology applied to Eq.(3), establishing the conditions under which the equation $I_{T}(x,T(y,z))=T(y,I_T(x,z))$ holds  holds becomes essential. Thus, we have the following results for a continuous Archimedean t-norm.
\\{\bf Proposition 5.10}  Let $T$ be a strict t-norm. $I_{T}(x,T(y,z))=T(y,I_T(x,z))$ holds for any $x>z$.\vspace{1mm}
\\{\bf Proof.}  As $x>z$ and $T$ is strict, we have $T(y,I_T(x,z))=t^{-1}\left(t(y)t\left(t^{-1}\left(\frac{t(z)}{t(x)}\right)\right)\right)=
t^{-1}\left(\frac{t(t^{-1}(t(y)t(z)))}{t(x)}\right)=I_{T}(x,T(y,z))$.
\\{\bf Proposition 5.11} Let $T$ be a nilpotent t-norm. If neither $y$ nor $z$ are zero divisors of $T$, then $I_{T}(x,T(y,z))=T(y,I_T(x,z))$ holds for any $x>z$.
\\{\bf Proof.} Since neither $y$ nor $z$ are zero divisors of the nilpotent t-norm $T$, we can obtain $I_{T}(x,T(y,z))$ $=t^{-1}((1-t(x)+t(t^{-1}(T_L(t(y),t(z)))))\wedge 1)=t^{-1}((1-t(x)+t(y)+t(z)-1)\wedge 1)=t^{-1}((-t(x)+t(y)+t(z))\wedge 1)=t^{-1}(-t(x)+t(y)+t(z))$. On the other hand, $T(y,I_{T}(x,z))=t^{-1}((t(y)+t(I_{T}(x,z))-1)\vee 0)=t^{-1}((t(y)+t(t^{-1}((1-t(x)+t(z))\wedge 1))-1)\vee 0)=t^{-1}(t(y)-t(x)+t(z))=I_{T}(x,T(y,z))$.
\\{\bf Theorem 5.12} Let a continuous Archimedean t-norm $T$ be used to interpret the logical connective ``AND" in the if-then rule. If $S_F(T(A_1, A_2)),T(A_1',A_2'))=T(S_F(A_1,A_1'),S_F(A_2,A_2'))$ holds, then $B'_1$ in Eq.(1) is $B'_1(y)=\sup\limits_{x_1\in U_1}I_T(S_F(A'_1,A_1),T(A_{1}(x_1),B'_{11}(y)))$ with $B'_{11}(y)=$\vspace{1mm} $\sup\limits_{x_2\in U_2}I_T(S_F(A'_{2},A_2),T(A_{2}(x_2),B(y))))$.\vspace{1mm}
 \\{\bf Proof.}  Since $I_T$ satisfies EP and $I_T(T(x,y),z)=I_T(x,I_T(y,z))$ holds for the t-norm $T$ according to Ref.\cite{Baczynski}, we have  $I_T(T(S_F(A_1',A_1),S_F(A_2',A_2)),T(T(A_{1}(x_1),A_{2}(x_2)),B(y)))=I_T(S_F(A_1',A_1),I_T$ $(S_F(A_2',A_2),T(A_{1}(x_1),I_T(A_{2}(x_2),B(y)))))$. By Propositions 5.10 and 5.11, it can be rewritten as $I_T(S_F(A_1',A_1),T(A_{1}(x_1),I_T(S_F(A_2',A_2),T(A_{2}(x_2),B(y)))))$ when $T(T$ $(A_{1}(x_1),A_{2}(x_2)),B(y))\neq 0$. Let $E_y=\{(x_1,x_2)|T(T(A_{1}(x_1),A_{2}(x_2)),B(y))\neq0\}$. Then, we have $B'_1(y)=\sup\limits_{(x_1, x_2)\in U_1\times U_2}I_T(T(S_F(A'_{1},A_1),S_F(A'_2,A_2)),T(T(A_{1}(x_1),A_{2}(x_2)),B(y)))=$\vspace{1mm} $\sup\limits_{(x_1, x_2)\in E_y}I_T(T(S_F(A'_{1},A_1),S_F(A'_2,A_2)),T(T(A_{1}(x_1),A_{2}(x_2)),B(y)))=\sup\limits_{(x_1, x_2)\in E_y}I_T(S_F(A_1',A_1),$\vspace{1mm} $T(A_{1}(x_1),I_T(S_F(A_2',A_2),T(A_{2}(x_2),B(y)))))=\sup\limits_{x_1\in U_1}I_T(S_F(A_1',A_1),T(A_{1}(x_1),\sup\limits_{x_2\in U_2}I_T(S_F(A_2',$\vspace{1mm} $A_2),T(A_{2}(x_2),B(y)))))$. Let  $B'_{11}(y)=\sup\limits_{x_2\in U_2}I_T(S_F(A'_{2},A_2),T(A_{2}(x_2),B(y))))$. We can obtain\vspace{1mm} $B'_1(y)=\sup\limits_{x_1\in U_1}I_T(S_F(A'_1,A_1),I_T(A_{1}(x_1),B'_{11}(y)))$.\vspace{1mm}

We can also extend this hierarchical Raha's SBAR method  to $n$-dimensional input fuzzy system. An algorithm for this approach may be similarly formulated as follows.
\\{\bf Algorithm 2} Hierarchical Raha's SBAR method based on Eq.(1).

\textbf{Input:} $A_k'$, $A_k(k=1,2\cdots,n)$, $B$

\textbf{Output:} $B'$
\\ 1. \textbf{for} $i=1$ to $n$ and $j=2$ to $n$ \textbf{do}
\\  2.  \quad\ \ $s_{i}=S_F(A_i',A_i)$
\\ 3. \quad\ \ $B'_{j}(y)=\inf\limits_{x_{i}\in U_{i}}I_T(s_{i},T(A_{i}(x_{i}),B(y)))$\vspace{1mm}
\\ 4. \quad\ \ $B\leftarrow B'_{j}$
\\ 5. \textbf{end for}
\\ 6. \textbf{return} $B'=B$

The computational and space complexity of Algorithm 2 can be analyzed analogously, though a detailed examination is unnecessary here.

To give a good view of the hierarchical Raha's SBAR method, we will respectively calculate the output with the hierarchical Raha's SBAR method and Raha's SBAR method in the following example.
\\{\bf Example 5.13}  Let $A_1=[1,0.9,0.6,0.7]$, $A_2=[0.4,0.4,0.6,0.5,0.3]$, $A_3=[0.6,0.3,0.5]$ and $B=[0.3,0.4,0.2,0.1]$. Suppose that the logical connective ``AND" in the if-then rule is interpreted by $T_P$ and $I=I_{GG}$ in Eq.(2).

For a fuzzy input $A'=A'_1\times A'_2\times A'_3$ with $A'_1=[0.8,0.5,0.7,0.9]$, $A'_2=[0.5,0.6,0.7,0.4,0.4]$ and $A'_3=[0.8,0.7,0.9]$, with  Raha's SBAR method, the output $B'_2$ is computed according to the following steps.

Step 1. It needs to compute $S_F(A'_1,A_1)=\frac{5}{9}$, $S_F(A'_2,A_2)=\frac{2}{3}$, $S_F(A'_3,A_3)=\frac{3}{7}$ and $s=T_P(S_F(A'_1,A_1),T_P(S_F(A'_2,A_2),S_F(A'_3,A_3)))=\frac{10}{63}$.\vspace{1.5mm}

Step 2. Calculate $T_P(A_1(x_1),A_2(x_2))=\left(
                 \begin{array}{ccccc}
                   0.4& 0.4& 0.6& 0.5&0.3 \\
                   0.36& 0.36& 0.54& 0.45&0.27 \\
                  0.24& 0.24& 0.3& 0.25&0.15\\
                   0.28& 0.28& 0.42& 0.35&0.21\\
                 \end{array}
               \right)$, $a=\mathop{\sup}\limits_{\mathbf{x}\in U}T_P$\vspace{1.5mm} $(T_P(A_1(x_1),A_2(x_2)),A_3(x_3))=0.324$ and $I_{GG}(a,B(y))=[\frac{5}{9}, 1, \frac{50}{81}, \frac{25}{81}]$.

Step 3.  Obtain $B'_2(y)=I_{GG}(s,I_{GG}(a,B(y)))=[1,1,1,1]$.  The computational count of $B'_2$ with Raha's SBAR method can be listed in Table 3.
$$\mbox{\bf{\small Table\ 3 \ The computational count of Raha's SBAR method}}$$
\begin{center}
\tabcolsep 0.05in
\begin{tabular}{cc}
 \toprule[1pt]
Procedure & Calculation count\vspace{1mm}\\
\midrule[0.75pt]
  $S_F(A_1',A_1)$       &$3\times4+3=15$\vspace{1mm}\\
  $S_F(A_2',A_2)$   &$3\times 5+4=19$\vspace{1mm}\\
  $S_F(A_3',A_3)$   &$3  \times 3+2=11$\vspace{1mm}\\
  $s$       &2\vspace{1mm}\\
  $T_P(A_1,A_2)$ &$4 \times 5=20$\vspace{1mm}\\
  $T_P(T_P(A_1,A_2),A_3)$ &$20 \times 3=60$\vspace{1mm}\\
  $a$ &59\vspace{1mm}\\
  $I_{GG}(a,B)$ &4\vspace{1mm}\\
 $B'_2$&4\vspace{1mm}\\
 Total &194\vspace{1mm}\\
   \bottomrule[1pt]
\end{tabular}
\end{center}\vspace{3mm}

Next, we compute the output $B'_2$ with our proposed hierarchical Raha's SBAR method.  Firstly, compute $S_F(A'_1,A_1)=\frac{5}{9}$, $S_F(A'_2,A_2)=\frac{2}{3}$ and $S_F(A'_3,A_3)=\frac{3}{7}$. Secondly, calculate $\mathop{\sup}\limits_{x_1\in U_1}A_1(x_1)=1$,  $\mathop{\sup}\limits_{x_2\in U_2}A_2(x_2)=$\vspace{1mm} $0.6$ and
 $\mathop{\sup}\limits_{x_3\in U_3}A_3(x_3)=0.6$. Further, we compute $B'_{23}(y)=\mathop{\inf}\limits_{x_3 \in U_3}I_{GG}(S_F(A'_3,A_3),I_{GG}(A_3(x_3),$\vspace{1mm} $B(y)))=[1,1,1,1]$, $B'_{22}(y)=\mathop{\inf}\limits_{x_3 \in U_3}I_{GG}(S_F(A'_2,A_2),I_{GG}$ $(A_2(x_2),B'_{23}(y)))=[1,1,1,1]$
 and $B'_{21}(y)=\mathop{\inf}\limits_{x_1\in U_1}I_{GG}(S_F(A'_1,A_1),I_{GG}(A_1(x_1),B'_{22}(y)))=$\vspace{1mm} $[1,1,1,1]$. Finally, we have $B_2'=B'_{21}=[1,1,1,1]$. Similarly, the computational count of $B'_2$ with our proposed hierarchical Raha's SBAR method can be illustrated in Table 4.
 $$\mbox{\bf{\small Table\ 4 \ The computational count of hierarchical Raha's SBAR method}}$$
\begin{center}
\tabcolsep 0.05in
\begin{tabular}{cc}
 \toprule[1pt]
Procedure & Calculation count\vspace{1mm}\\
\midrule[0.75pt]
  $S_F(A_1',A_1)$       &$3\times4+3=15$\vspace{1mm}\\
  $S_F(A_2',A_2)$   &$3\times 5+4=19$\vspace{1mm}\\
  $S_F(A_3',A_3)$   &$3  \times 3+2=11$\vspace{1mm}\\
  $\mathop{\sup}\limits_{x_1\in U_1}A_1(x_1)$&3\vspace{1mm}\\
    $\mathop{\sup}\limits_{x_2\in U_2}A_2(x_2)$&4\vspace{1mm}\\
    $\mathop{\sup}\limits_{x_3\in U_3}A_3(x_3)$ &2\vspace{1mm}\\
  $B_{23}'$ &$4+4=8$\vspace{1mm}\\
   $B_{22}'$ &$4+4=8$\vspace{1mm}\\
 $B'_2$&4+4=8\vspace{1mm}\\
 Total &68\vspace{1mm}\\
   \bottomrule[1pt]
\end{tabular}
\end{center}\vspace{1mm}

Clearly, our proposed hierarchical Raha's SBAR method can effectively cut down the computational complexity. Moreover, it is enough to store $4+5+3$ memory cells for the antecedent and $4$ memory cells for the consequent of the if-then rule. And then we can quickly obtain the output $B_2'$ according to Algorithm 1.
\section{Conclusions}
\qquad
To facilitate the application of  REF in approximate reasoning, this work has studied the performance  of Raha's SBAR method employing REF. The conclusions mainly include

(1)  Characterizing REF by a given
binary aggregation function.

(2) Constructing REF with a given aggregation function and a mapping $f:[0,1]^2\rightarrow [0,1]$ satisfying
I1, CC and OP.

(3) Defining the approximate equality of two fuzzy sets with REF, and subsequently investigating the approximation equality of Raha's SBAR method with REF.

(4) Proposing two hierarchical Raha's SBAR methods with REF.

(5) Comparing the computational and space complexity of Raha's SBAR method with that of our proposed hierarchical Raha's SBAR methods.
\section{Directions for future research}
\qquad Our investigations can be applied in the following problems
\begin{itemize}
  \item Constructing REF with some aggregation functions.
    \item Comparing two images with REF.
  \item Classifying some high-dimensional images using the hierarchical Raha's  SBAR method.
\end{itemize}
These results show some application of REF in approximate reasoning. To be more practical, we could investigate the impact different REFs on the approximation
accuracy of fussy system using the hierarchical Raha¡¯s SBAR method.
\section{Acknowledgements}\qquad This work was funded by the National Natural
Science Foundation of China (Grant No. 61673352).

\end{document}